%% file: main.tex
\documentclass[review,onefignum,onetabnum]{siamart190516}
\nolinenumbers

\usepackage[colorinlistoftodos,bordercolor=orange,backgroundcolor=orange!20,linecolor=orange,textsize=scriptsize]{todonotes}

\usepackage{graphicx}
\usepackage{float}
\usepackage{subfigure}
\usepackage{tabularx}
\usepackage{multirow}

\usepackage{latexsym}

\usepackage{amsmath,amstext,amssymb,amsbsy}
\usepackage{mathtools}
\usepackage{bm,bbm}
\usepackage{caption}
\usepackage{subcaption}
\usepackage{xcolor}

\usepackage{cleveref}

\usepackage{algorithm}  
\usepackage{algorithmicx}  
\usepackage{algpseudocode}

\usepackage{alphalph}

\usepackage[multiple]{footmisc}

\newcommand{\Ll}{\mathcal{L}}
\newcommand{\B}{\mathcal{B}}
\newcommand{\K}{\mathcal{K}}

\newcommand{\Am}{\mathbf{A}}
\newcommand{\Mm}{\mathbf{M}}
\newcommand{\Rm}{\mathbf{R}}
\newcommand{\Pm}{\mathbf{P}}
\newcommand{\fv}{\mathbf{f}}
\newcommand{\uv}{\mathbf{u}}
\newcommand{\rv}{\mathbf{r}}
\newcommand{\zv}{\mathbf{z}}
\newcommand{\xv}{\mathbf{x}}

\newcommand{\NO}{NeuralOp}
\newcommand{\CO}{Class\_ite}
\newcommand{\kr}{Krylov\_ite}
\newcommand{\R}{\mathbb{R}}
\newcommand{\Pp}{\mathcal{P}}
\newcommand{\F}{\mathcal{F}}
\newcommand{\M}{\mathcal{M}}
\newcommand{\qqed}{\hfill\blacksquare}

\newcommand{\Omegab}{\bar{\Omega}}

\newtheorem{pref}{Preference function} %Test Case}
\newtheorem{remark}{Remark}
\newtheorem{assumption}{Assumption}[section]
\newtheoremstyle{assumption}
  {6pt}% measure of space to leave above the theorem. E.g.: 3pt
  {6pt}% measure of space to leave below the theorem. E.g.: 3pt
  {\rm}% name of font to use in the body of the theorem
  {}% measure of space to indent
  {\bfseries}% name of head font
  {}% punctuation between head and body
  {}% Manually specify head
\theoremstyle{assumption}

\DeclareMathOperator*{\argmin}{arg\,min}

\title{Automatic discovery of optimal meta-solvers via multi-objective optimization
\thanks{Submitted to the editors \today. Youngkyu Lee and Shanqing Liu contributed equally to this work. 
\funding{This work was supported by DARPA-DIAL grant HR00112490484.
.}}}

\author{Youngkyu Lee\footnotemark[2]  \ \and
Shanqing Liu\footnotemark[2] \ \and
J\'er\^ome Darbon\footnotemark[2] \footnotemark[3]  \and
George Em Karniadakis\footnotemark[2]\ %\footnotemark[3]
}

\begin{document}

\maketitle
\renewcommand{\thefootnote}{\fnsymbol{footnote}}
\footnotetext[2]{Division of Applied Mathematics, Brown University (\{youngkyu\_lee,  shanqing\_liu, jerome\_darbon, george\_karniadakis\}@brown.edu).}
\footnotetext[3]{Corresponding author.}

\begin{abstract}
We design two classes of ultra-fast meta-solvers for linear systems arising after discretizing PDEs by combining neural operators with either simple iterative solvers, e.g., Jacobi and Gauss-Seidel, or with Krylov methods, e.g., GMRES and BiCGStab, using the trunk basis of DeepONet as a coarse preconditioner. The idea is to leverage the spectral bias of neural networks to account for the lower part of the spectrum in the error distribution while the upper part is handled easily and inexpensively using relaxation methods or  fine-scale preconditioners.
We create a pareto front of optimal meta-solvers using a plurarilty of metrics, and we introduce a preference function to select the best solver most suitable for a specific scenario.
This automation for finding optimal solvers can be extended to nonlinear systems and other setups, e.g. finding the best meta-solver for space-time in time-dependent PDEs.
\end{abstract}



\input{introduction}
\input{preliminary}
\input{methodology}

\input{numerical_poisson}

\input{conclusion}

\newcommand{\etalchar}[1]{$^{#1}$}

\input{appendix}

\end{document}

%% file: introduction.tex
\section{Introduction}
\label{sec:introduction}
\subsection{Motivation and Background}

Partial differential equations (PDEs) form a cornerstone of mathematical modeling in a wide range of scientific and engineering disciplines. 
%Numerical solution of PDEs involves discretization techniques which leads to solving linear or nonlinear finite-dimensional system of equations. 
A long-standing challenge in scientific computing is to find fast and accurate numerical methods for numerically solving PDEs, especially at scale when billions of unknowns are involved. An effective solution is even more crucial in multi-query applications, e.g., 
uncertainty quantification, where the PDE needs to be solved many times for different initial and boundary conditions, different material properties or different excitation terms. 

Numerical schemes for solving PDEs have been indispensable tools in computational science and engineering since their pioneering development more than half a century ago. Among the most commonly used numerical techniques are finite difference methods~\cite{strikwerda2004finite}, finite element methods~\cite{bathe2006finite}, and spectral/spectral-element  methods~\cite{karniadakis2005spectral}. Upon any type of discretization,  a linear system of discrete equations needs to be solved. Before the 1970s, these linear systems used to be solved using simple iterative solvers, such as Jacobi, Gauss-Seidel methods (see, for instance,~\cite{greenbaum1997iterative,saad2003iterative}) and their variants Successive Over-Relaxation (SOR) and Symmetric SOR (SSOR) (\cite{berman1994nonnegative, hadjidimos2000successive}) methods. 
Today, for large-scale linear systems, Krylov methods, such as Conjugate Gradient (CG), Generalized Minimal Residual (GMRES) and Biconjugate gradient stabilized (BiCGStab) methods are more commonly employed (see, for instance,~\cite{saad1993flexible,van2003iterative, van1992bi}). These methods are favored for their lower computational cost per iteration and their suitability for parallel computing environments. 
Advanced techniques have been proposed, such as multi-grid method (see for instance~\cite{bramble1990parallel,hackbusch2013multi,wesseling2004introduction}) and domain decomposition methods (see for instance~\cite{mathew2008domain,xu1992iterative}).

Recently, scientific machine learning (SciML) has gained significant attention in numerically solving PDEs arising from physics and engineering. 
Since the groundbreaking work in~\cite{raissi2019physics}, various frameworks have been proposed, among which we mention the deep operator network (DeepONet)~\cite{lu2021learning}, U-DeepONet~\cite{diab2024u}, Fourier Neural Operator (FNO)~\cite{li2020fourier}, Deep operator network based on Kolmogorov Arnold networks (DeepOKAN)~\cite{abueidda2024deepokan}, and Transformer-based neural operators~\cite{shih2024transformers}. 
PINNs and neural operators have shown promise in diverse fields, including but not limited to, fluid mechanics~\cite{raissi2020hidden, cai2021physics}, material science~\cite{zhang2022analyses, goswami2022physics}, electro dynamics~\cite{cai2021deepm} and optimal controls~\cite{meng2022sympocnet, zhang2024time}. 

For the neural operator approaches, approximating the operator of the discretized system at the fine level is challenging due to the spectral bias of neural networks . To address this difficulty, hybrid approaches have been introduced in~\cite{zhang2024blending, kahana2023geometry, kopanivcakova2024deeponet}, where neural operators are used to define preconditioners for iterative solvers, rather than serving as direct solvers themselves. These hybrid approaches have proven effective in reducing computational time and the number of iterations.

Using the hybrid preconditioning strategies, several families of meta-solvers can be generated. %: one blending a relaxation iterative solver with a neural operator, and the other combining a Krylov method with a neural operator. 
This raises a natural question: which meta-solver is the ``optimal'' one? 
From the user's perspective, the answer depends on the specific task and context, as various factors such as speed, accuracy, memory usage, and other performance dimensions must be considered when evaluating different solvers. 
The main objective of this paper is to address this question by proposing a methodology to automatically discover the optimal meta-solver for a given target problem, via tools from multi-objective optimization~\cite{miettinen1999nonlinear}. 

\subsection{Contribution}
In this paper, we introduced a parameterization method for the meta-solvers within subspace $\M^r$ and $\M^K$, for  relaxation-based and Krylov-based methods, respectively. A meta-solver is instantiated by giving a specific element from this subspace. 
The performance of the meta-solvers is then quantified using vector-valued maps from $\M^r$ and $\M^K$ to $\R^N$, which captures and records the performance criteria, such as computational time, relative error, and memory allocations. Two families of meta-solvers are then generated and tested for solving 1-d, 2-d and 3-d Poisson equations. %Most performance criteria provide significant variability that allows us to discriminate between meta-solvers. 

We developed a mathematical formulation to identify the ``optimal'' solvers. 
In particular, the evaluation of the performance of different solvers is formulated as a multi-objective optimization (MOO) problem~\cite{{miettinen1999nonlinear}}, in the sense of minimizing the vector-valued objective maps. We defined the concepts of dominance, Pareto optimality~\cite{censor1977pareto}, and Pareto front for the meta-solvers. 
We developed and implemented a MOO solver to numerically identify Pareto optimal solvers. 

We also developed a novel preference function-based discovery methodology for identifying optimal solvers among Pareto ones. Once users provide their ``preferences'', the discovery of the optimal solver is reformulated as a classical optimization problem: $\min p \circ f$ over the Pareto optimal solvers, where $p: \R^N \to \R$ models the preference of a user and $f$ is the vector-valued objective map. We present the numerical results for a particular family of preference function: the weighted sum function. 
Moreover, we develop a Linear Programming (LP)~\cite{dantzig2002linear} based algorithm for re-discovering a particular known meta-solver among the Pareto ones. %, that is, g
Given a meta-solver $x$, the algorithm identifies the weights that makes $x$ the optimal one. 

This paper is organized as follows. \Cref{sec:preliminaries} provides an overview of the background information, including problem formulation in \Cref{subsec::problem}, classical iterative solvers in \Cref{subsec::iterative}, and hybrid preconditioning using neural operators in \Cref{subsec::hybrid}. 
In~\Cref{sec-method}, we introduce our multi-objective optimization and Pareto optimality-based evaluation methodology for parameterized meta-solvers. 
We also present the approach for discovering optimal solvers using preference function in~\Cref{subsec-discovery}, and the re-discovery of optimal solvers via Linear Programming in~\Cref{sec-red}. 
Numerical results are presented in~\cref{sec-num}, with a focus on relaxation-based methods in \Cref{subsec-relax} and Krylov-based methods in \Cref{sec-poi-kry}. 
Finally, we summarize the findings and present concluding remarks in \Cref{sec-conclu}. 

%% file: preliminary.tex
\section{Preliminaries}
\label{sec:preliminaries}
This section provides some background materials used in the reminder of this paper.

\subsection{Problem Formulation}
\label{subsec::problem}
As a starting point of this study, we are interested in numerically approximating the solution $u : \Omegab \to \R$ of the linear differential equation of the form
\begin{equation}\label{eq:linear}
\left\{
\begin{aligned}
& \Ll(u) = f, \ &\text{in } \Omega \ , \\
& \B(u) = g, \ &\text{in } \partial \Omega \ ,
\end{aligned}
\right.
\end{equation}
where $\Omega \subset \R^d$ is an open bounded domain, $\Ll : \R^{\Omega} \to \R^{\Omega}$ is a linear differential operator.
$f:\Omega \to \R$ is a known function that does not explicitly involve $u$. $\B :\R^{\partial \Omega} \to \R^{\partial \Omega}$, $g:\partial \Omega \to \R$ serve as the boundary condition. 
Given a mesh $\Omega^h$ discretizing $\Omega$, a common approach to numerically approximate the solution of system~\eqref{eq:linear} is the 
%Problem~\eqref{eq:linear} is typically 
%using 
finite element method, which can be formulated as a linear system on $\Omega^h$, of the form
\begin{equation}
\label{eq:fem}
    \Am \uv = \fv \ ,
\end{equation}
where $\uv$ denotes the nodal coefficients of finite element basis functions. Note that if the PDE is nonlinear, then linearization techniques can be used to reduce the problem in this form during the 
solution of the nonlinear algebraic system or during time-stepping e.g. employing implicit-explicit (IMEX) methods.

\subsection{Iterative Solvers and Preconditioners}
\label{subsec::iterative}
The discretized system~\eqref{eq:fem} can be solved using iterative solvers, e.g., relaxation methods and Krylov methods~\cite{saad2003iterative}. 
In particular, given an appropriate initialization $\uv^{(0)}$, 
let $\uv^{(i)}$ be the approximate solution after $i$ iterations.  
At the $(i+1)$-th step, the interation has the form 
\begin{equation}\label{eq::ite}
\left\{ 
\begin{aligned}
& \rv^{(i)} = \fv - \Am \uv^{(i)} \ , \\
& \uv^{(i+1)} = \uv^{(i)} +  \K (\rv^{(i)}) \ ,
\end{aligned}
\right.
\end{equation}
where $\K$ denotes the update process, which depends on the iterative solver. 
Classical iterative solvers can be accelerated by utilizing certain preconditioners. 
At the $(i+1)$-th iterate, the preconditioning process is conducted by
\begin{equation}
\label{eq:preconditioner}
\left\{
    \begin{aligned}
        &\rv^{(i)} = \fv - \Am \uv^{(i)} \ , \\
        &\zv^{(i)} = \Mm(\rv^{(i)}) \ , \\
        &\uv^{(i+1)} = \uv^{(i)} + \K (\zv^{(i)}) \ ,
    \end{aligned}
\right.
\end{equation}
where $\Mm$ denotes the preconditioner. 
Since the convergence rate of iterative solver is dependent on the condition number $\kappa$ of given system, the purpose of the preconditioner is $\kappa (\Mm \Am) \approx 1$ or $\kappa (\Am \Mm) \approx 1$.
Typically, the incomplete LU decomposition~(ILU) is used as the preconditioner when there is no prior information about the domain and equation.

\subsection{Neural Operators and Hybrid Preconditioners}
\label{subsec::hybrid}
%On the other hand, t
Another approach to solving problem~\eqref{eq:linear} is through the use of neural operators~\cite{kovachki2021neural,lu2021learning}, which approximate the functional operator that maps function $f$ to the solution $u$ using neural networks. 
However, approximating this operator at the fine level is challenging due to the spectral bias of neural networks~\cite{rahaman2019spectral}. 
To address this difficulty, 
hybrid approaches have been introduced in~\cite{kopanivcakova2024deeponet, zhang2024blending}, where neural operators are used to define preconditioners for iterative solvers, rather than using neural operators as direct solvers for~\eqref{eq:linear}. 
%hybrid approaches that use the neural operators to define the preconditioners of iterative solvers are introduced in~\cite{kopanivcakova2024deeponet,zhang2024blending}, 
%instead of using neural operator as a direct solver of~\eqref{eq:linear}. %, hybrid approaches that use the neural operators to define the preconditioners of iterative solvers are introduced in~\cite{kopanivcakova2024deeponet,zhang2024blending}.
In this context, 
the preconditioning process~\eqref{eq:preconditioner} can be further accelerated as
\begin{equation}
\left\{
    \begin{aligned}
        &\rv^{(i)} = \fv - \Am \uv^{(i)} \ , \\
        &\uv^{(i+1/2)} = \uv^{(i)} + \Mm_{1}(\rv^{(i)}) \ , \\
        &\rv^{(i+1/2)} = \fv - \Am \uv^{(i+1/2)} \ , \\
        &\zv^{(i+1/2)} = \Mm_{2}(\rv^{(i+1/2)}) \ , \\
        &\uv^{(i+1)}   = \uv^{(i+1/2)}+ \K (\zv^{(i+1/2)}) \ ,
    \end{aligned}
\right.
\label{numerical_scheme}
\end{equation}
where $\Mm_{1}$ and $\Mm_{2}$ denote the steps of relaxation method and an inference through the pre-trained neural operator, respectively. 

System~\eqref{numerical_scheme} is typically well defined 
for relaxation methods based hybrid preconditioners as in~\cite{zhang2024blending}. 
However, when the preconditioner $\Mm_{2}$ is defined as the inference via the neural operator, it cannot be directly used as the preconditioner for Krylov methods.
%In order to handle this problem, 
To address this issue, when DeepONet~\cite{lu2021learning} is utilized as the neural operator, the trunk-basis~(TB) approach, proposed in~\cite{kopanivcakova2024deeponet}, extracts the prolongation operator~$\Pm$ and restriction operator~$\Rm$ from the trunk network of DeepONet. 
This leads to 
the construction of second preconditioner $\Mm_{2}:=\Pm(\Rm\Am\Pm)^{-1}\Rm=\Pm\Am_{c}^{-1}\Rm$.
The $(i,j)$-th component of the matrix of the prolongation operator~$\Pm$ is computed by
\begin{equation}
    [\Pm]_{ij} = T_{j}(\xv_{i}) \ ,
\end{equation}
where $T_{j}(\xv_{i})$ denotes the $j$-th component of the output of the trunk network evaluated at the coordinate $\xv_{i}$.
Note that the restriction operator is defined as the adjoint operator of prolongation operator, i.e., $\Rm = \Pm^{\ast}$.
The quality of the prolongation operator $\Pm$ can be further improved utilizing the sampling strategies and QR decomposition (see~\cite{kopanivcakova2024deeponet}).

%% file: methodology.tex
\section{Multi-objective Optimization (MOO), Pareto Optimal, Discovery and Re-discovery of Optimal Solvers}
\label{sec-method}

In this section, we present our approach for identifying optimal meta-solvers using tools from multi-objective optimization. 
For the meta-solvers, we consider both the relaxation-based methods and the Krylov-based methods, as presented in~\Cref{subsec::hybrid}. 
In particular, we present the discovery and re-discovery of optimal solvers based on the Pareto optimality and preference functions. 

\subsection{Parameterization of Meta-Solvers}
\begin{figure}
\centering
\subfigure[Relaxation-based meta-solvers]{
\includegraphics[width=0.45\textwidth]{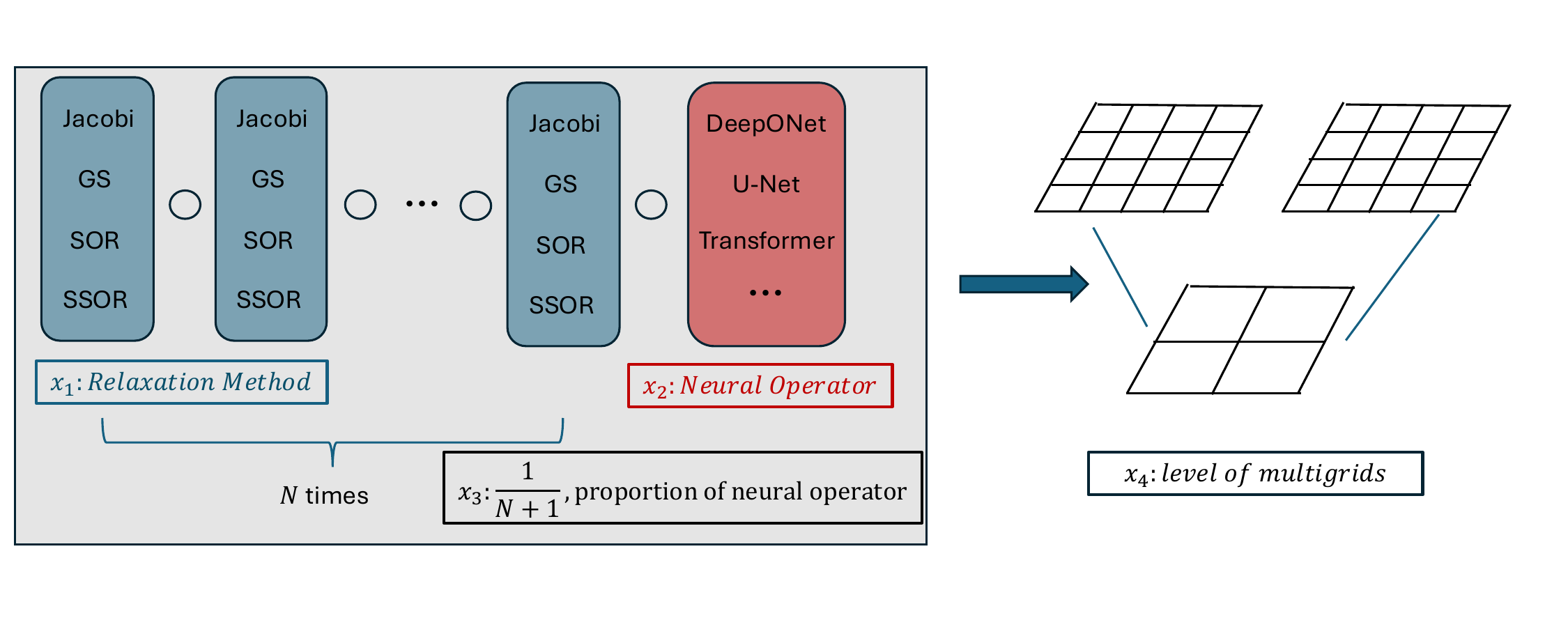}}
\subfigure[Krylov-based meta-solvers]{
\includegraphics[width=0.45\textwidth]{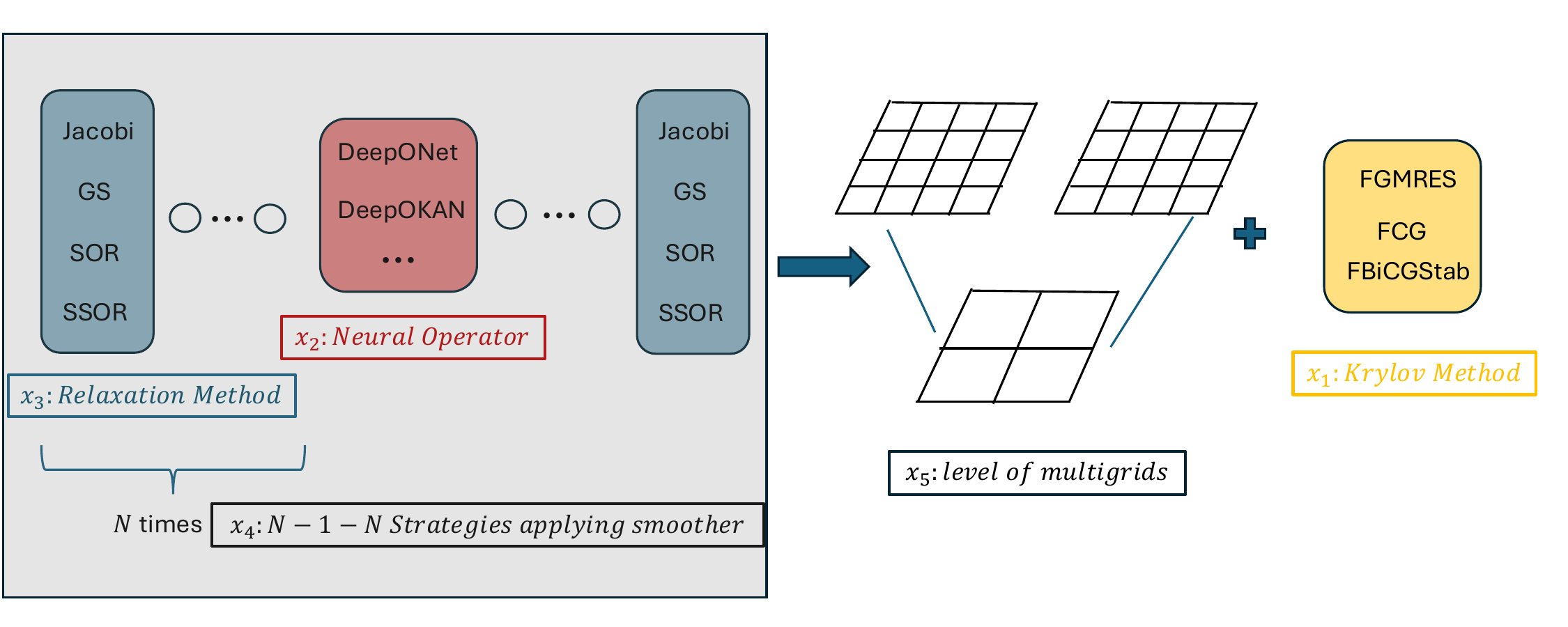}}
\caption{Construction of meta-solvers. A relaxation-based meta-solver combines a neural operator with an iterative solver, instantiated using a fixed proportion. Additionally, we apply a multi-grid technique on top of this combination. For the Krylov-based meta-solver, we first select a relaxation method as a smoother and combine it with a neural operator using the strategy: $N$ steps of smoother, followed by 1 step of the neural operator, and then $N$ more steps of the smoother. Additionally, we apply a multi-grid technique. Finally, a Krylov method is applied on top of this framework.}
\label{fig::sketch_construc}
\end{figure}

We begin by defining the parameter space of meta-solvers. 
A sketch of the construction of relaxation-based methods, and of Krylov-based methods are given in~\Cref{fig::sketch_construc}.
\subsubsection{Relaxation-based Methods}
In the context of relaxation-based methods, a meta-solver is a combination of different neural operators with classical iterative solvers, which is eventually instantiated using different proportion between the two components. We also explore the creation of meta-solver using some advanced techniques, starting with the multigrid method. 
The meta-solvers are then parameterized within a four-dimensional subspace $\M^r$, 
where each dimension represents the choice of neural operators, the choice of classical iterative solvers, proportion of neural operator, and levels of multi-grid techniques, respectively.  

For neural operators, we consider DeepONet, U-Net DeepONet, Fourier Neural Operator (FNO),  Transformer, Kolmogorov-Arnold Networks (KAN), Jacobi KAN and Chebyshev KAN (See~\cite{lu2021learning} for DeepONet, \cite{diab2024u} for U-Net, \cite{li2020fourier} for FNO, \cite{shih2024transformers} for Transformer, and~\cite{shukla2024comprehensive} for the family of KANs). As for classical numerical iterators, we consider standard relaxation methods, such as Jacobi method, Gauss-Seidel method, Successive Over-Relaxation (SOR) method and Symmetric Successive Over-Relaxation (SSOR). The proportion of the neural operator to classical iterative operator ranges from $\frac{1}{128}$ to $\frac{1}{2}$. Among all versions of the multi-grid technique, we first consider 1-level (Non multi-grid), 2-level and 3-level multi-grids.

To summarize, we parameterize the relaxation-based meta-solvers by the variable $x =  (x_1,x_2,x_3,x_4)$ such that 
\begin{equation}\label{relax_para_solver}
	\left\{
	\begin{aligned}
		& x_1 \in \NO = \left\{ \text{DeepONet, U-Net, FNO, Transformer,} \right. \\
   &\qquad \qquad \qquad\qquad \qquad \qquad \left. \text{ KAN, Jacobi KAN, Cheby KAN} \right\}   ;  \\
		& x_2 \in \CO = \{ \text{Jacobi, Gauss-Seidel, SOR, SSOR}    \} ; \\ 
            & x_3 \in P = \left\{\frac{1}{128} , \frac{1}{64}, \frac{1}{32}, \frac{1}{16}, \frac{1}{8},\frac{1}{4},\frac{1}{2}\right\}, \ \text{the proportion of neural operator}  ; \\
            & x_4 \in \{0,1,2\}, \ \text{levels of multi-grid}. 
           % & x_4 \in \left\{\frac{1}{16}, \frac{1}{32}, \frac{1}{64}, \frac{1}{128} \right\}, \ \text{mesh step}; \\
		\end{aligned}
	 \right.
	\end{equation}
Therefore, a meta-solver is instantiated by giving a specific $x \in \M^r:= \NO \times \CO  \times P \times \{0,1,2\}$.
\subsubsection{Krylov-based Methods}
The Krylov-based methods are considered in the trunk basis hybridization approaches, that are combinations of different neural operators and Krylov methods. On top of the neural operators, we apply several steps of relaxation methods as smoother. 
In more details, considering the abstract numerical schemes~\eqref{numerical_scheme}, $\Mm_1$ is the relaxation method as smoother and 
$\Mm_2$ is the Krylov-based preconditioner constructed by trunk basis hybridization approaches.  
%and $\Mm_1$ is the relaxation method as smoother.  
We adapt various strategies of applying smoothers. 
Moreover, the creation of Krylov-based methods also consist of adapting multigrid method. Therefore, in this case the meta-solvers are parameterized within a five-dimensional subspace $\M^K$, where each dimension represents the choice of neural operators, the choice of Krylov methods, the choice of relaxation methods as smoothers, strategies of applying the smoothers, and the levels of multi-grid techniques, respectively. 

For neural operators, we consider DeepONet, U-Net DeepONet, Kolmogorov-Arnold Networks (KAN), Jacobi KAN and Chebyshev KAN. We consider three types of Krylov method: flexible generalized minimum residual (FGMRES), flexible conjugate gradients (FCG) and flexible biconjugate gradient stabilized (FBiCGStab) iterative methods. For relaxation method, we consider Jacobi method, Gauss-Seidel method, Successive Over-Relaxation (SOR) method and Symmetric Successive Over-Relaxation (SSOR) as smoothers. The strategies of applying smoothers range from $1-1-1$ to $9-1-9$. The levels of multi-grid are choosen from $1$, $2$ and $3$. %A sketch of the Krylov-based methods is given in~\Cref{fig::krylov}.
%\begin{figure}
%\centering
%\includegraphics[width=0.8\textwidth]{figs/sketch_krylov.pdf}
%\caption{Krylov-based meta-solvers}
%\label{fig::krylov}
%\end{figure}
To summarize, we parameterize the Krylov-based meta-solvers by the variable $x =  (x_1,x_2,x_3,x_4,x_5)$ such that 
\begin{equation}\label{krylov_para_solver}
	\left\{
	\begin{aligned}
		& x_1 \in \NO = \left\{ \text{DeepONet, U-Net, KAN, Jacobi KAN, Cheby KAN} \right\}  ;  \\
            & x_2 \in \kr  = \{ \text{FGMRES, FCG, FBiCGStab} \}  ; \\ 
		& x_3 \in \CO = \{ \text{Jacobi, Gauss-Seidel, SOR, SSOR}    \} ; \\ 
            & x_4 \in S =  \{ \text{1-1-1, 3-1-3, 5-1-5, 7-1-7, 9-1-9 }\}, \ \text{strategies applying smoother} ; \\
            & x_5 \in \{0,1,2\}, \ \text{levels of multi-grid}. 
           % & x_4 \in \left\{\frac{1}{16}, \frac{1}{32}, \frac{1}{64}, \frac{1}{128} \right\}, \ \text{mesh step}; \\
		\end{aligned}
	 \right.
	\end{equation}
 A Krylov method based meta-solver is instantiated by giving a specific $x \in \M^K:= \NO \times \kr \times \CO \times S \times \{0,1,2\}$.

 \subsection{Performance Evaluation by Multi-objective Optimization (MOO)}
 Once a meta-solver is constructed, its performance is evaluated using several criteria. Herein, we consider seven criteria after a meta-solver has been applied to solve specific problems/equations: speed, measured in terms of total computational (execution) time; accuracy, evaluated in terms of relative error; convergence rate, computed in the limit sense, taking the inverse of the number; memory allocation, that is referring to the maximum memory required during the entire execution of the solver; multiply-accumulate-operations (MACs); and the training time of the neural network. Each criterion is quantified by a real number. 
 Considering the relaxation-based methods, given a meta-solver $x \in \M^r$, its performance is evaluated using a vector-valued map $f^r: \M^r \to \R^7$,  which maps $x$ to $(f_1(x), f_2(x),\dots, f_7(x))$ with $f_1(x)$ is the computational time, $f_2(x) $ is the relative error, $f_3(x)$ is the number of iterations, $f_4(x) $ is the convergence rate, $f_5(x) $ is the  memory allocation, $f_6(x)$ is the MACs and $f_7(x)$ is the training time, respectively.
%\begin{equation}
	%\left\{
	%\begin{aligned}
	%	& f_1(x) , \text{ computational time}; \\
	%	& f_2(x) , \text{ relative error}; \\ 
	%	& f_3(x) , \text{ number of iteration};  \\
	%	& f_4(x) , \text{ convergence rate}; \\ 
	%	& f_5(x) , \text{ memory allocation}; \\
        %   & f_6(x) , \text{ MACs}; \\
        %   & f_7(x) , \text{ training time}. 
	%	\end{aligned}
	%\right.
	%\end{equation}
    
The same performance criteria apply to the Krylov-based method, which is a map with the same image set as $f^r$, and is defined in $\M^K$, $f^K: \M^K \to \R^7$. Note that the methodology proposed here and below holds, {\em mutatis mutandis}, allowing for the addition of new criteria or the removal of existing ones among $f_1$ to $f_7$. 

Finding the optimal solver can be characterized as a multi-objective optimization problem, in the minimization sense, of the form
\begin{equation}\label{eq:moo}
 \min_x  \ f(x) = 
\left[ f_1(x) , \dots, f_7(x) \right]^T \ ,
\end{equation}
where for the relaxation-based method the minimization is taken over the domain $\M^r$, and for the Krylov-based method, the minimization is taken over the domain $\M^K$. 
However, problem~\eqref{eq:moo} does not necessarily admit a solution, that is, there is generally no $x$ that simultaneously minimizes all objective maps $f_1$ to $f_7$. This is consistent with the fact that different problems and criteria lead to different optimal solvers. This motivated us to further develop a tool to address this multi-objective optimization problem in order to identify the optimal solvers. 

\subsection{Identifying Optimal Solvers in Pareto Sense}

Consider a general multi\ -objective optimization problem, that is to minimize a $N-$dimensional vector valued function $L :\R^d \to \R^N$:
\begin{equation}\label{MOO_obj}
L(x) = \left[ 
L_1(x) 
\dots, 
L_N(x)
\right]^T
  \ ,
\end{equation}
over a domain $X \subseteq \R^d$. 
In general, the multi-objective optimization problem does not admit a solution, $x^*\in X$, that minimizes all the objective functions $\{L_i\}_{i \in \{1,2,\dots,N \}}$. Thus, a solution must be defined on certain concept of sub-optimality. Here, our objective is to identify candidates that are in the context of Pareto optimality.

We first recall the definition of dominance, in the sense of minimize $L$ over $X$.
\begin{definition}\label{MOO_dominance} For every $x^1, x^2 \in X$, 
    \begin{enumerate}
        \item we say that $x^1$ dominates $x^2$ if and only if $L_i(x^1) \leq L_i(x^2)$, for every $i \in \{1,2,\dots, N\}$, and there exists at least one $j$ such that $L_j(x^1)<L_j(x^2)$. We denote $x^1 \succeq x^2$ when $x^1$ dominates $x^2$. 
        \item we say that $x^1$ strictly dominates $x^2$ if and only if $L_i(x^1) < L_i(x^2)$, for every $i \in \{1,2,\dots, N\}$. We denote $x^1 \succ x^2$ when $x^1$ strictly dominates $x^2$. 
    \end{enumerate}
\end{definition}

The notion of dominant solution is closely related to Pareto optimality.
\begin{definition}\label{MOO_Pareto} For every $x\in X$,
\begin{enumerate}
    \item we say that $x$ is a strong Pareto optimal solution if and only if there is no solution in $X$ that dominates $x$. We denote $\Pp_L(X)$ the set of all strong Pareto optimal solution of minimizing $L$ over $X$. 
    \item we say that $x$ is a weak Pareto optimal solution if and only if there is no solution in $X$ that strictly dominates $x$. We denote $\Pp^w_L(X)$ the set of all weak Pareto optimal solution of minimizing $L$ over $X$.
\end{enumerate}
\end{definition}

It follows directly from~\Cref{MOO_Pareto} that $\Pp_L(X) \subseteq \Pp^w_L(X)$. 
\begin{definition}\label{MOO_ParetoFront} 
In the context of minimizing $L$ over $X$, we call the Pareto front $\F_L(X)$ (weak Pareto front $\F^w_L(X)$ resp.) the image of the Pareto optimal solution $\Pp_L(X)$ (weak Pareto optimal solution $\Pp^w_L(X)$ resp.) by the multi-objective map $L$, that is
    \begin{equation}\label{pareto_front}
    %\begin{equation}
    \F_L(X) := \{L(x) \mid x \in \Pp_L(X) \} , \  
    %\begin{equation}\label{weak_pareto_front}
    F^w_L(X) := \{L(x) \mid x \in \Pp^w_L(X) \} \ .
    %\end{equation}
    \end{equation}
\end{definition}
It is worth noticing that when there exists a unique solution $x^*$ that optimizes all the objective functions $\{L_i\}_{i \in \{1,2,\dots,N \}}$, the Pareto front is a singleton
\begin{equation}
\F_L(X) = \F^w_L(X) = \{L_1(x^*), L_2(x^*), \dots, L_N(x^*)\} \ .
\end{equation} 

The Pareto optimal solutions to problem~\eqref{eq:moo}  correspond to the optimal meta-solvers in the Pareto sense. %A sketch illustrating the identification of Pareto optimal solvers is shown in \Cref{fig::front}. 

Note that the Pareto optimality framework applies when one adjusts the parameterization of the solvers, or the performance matrix. In these cases, it is sufficient to modify the domain and image set of $f$. 
In particular cases, we consider the following assumption.
\begin{assumption}\label{assump1}
For every $x_i,x_j \in \Pp_L(X)$ and $x^i \neq x^j$, there exist at least one $k \in \{1,2,\dots,N\}$ such that $L_k(x^i) \neq L_k(x^j)$.
\end{assumption}

~\Cref{assump1} naturally holds in the context of evaluating different meta-solvers, as it is rare to expect that two different solvers would exhibit identical performance across all criteria. Under~\Cref{assump1}, when additional solvers or performance criteria are added in the existing ones, the following property of the set of Pareto optimal solvers holds. 

\begin{proposition}\label{mono_proper}
    \begin{enumerate}
    \item\label{mono_l} For any $\tilde{L}: \R^d \to \R^{N+\tilde{N}}$ such that $x \to \tilde{L}(x) = (L_1(x), \dots, L_N(x), \tilde{L}_{N+1}(x) \dots, \tilde{L}_{N + \tilde{N}})$,  %denote $\Pp_{\tilde{L}}(X)$ ($\Pp^w_{\tilde{L}}(X)$ resp.) the set of Pareto (weak Pareto resp.) optimal solutions in the sense of minimizing $\tilde{L}$ over $X$. W
    we have \begin{equation}
    \Pp_L(X) \subseteq \Pp_{\tilde{L}}(X), \quad \Pp^w_L(X) \subseteq \Pp^w_{\tilde{L}}(X)\ . 
    \end{equation}
    \item\label{mono_x} For any $\tilde{X} \subseteq \R^d$ such that $X \subseteq \tilde{X}$, %denote $\tilde{\Pp}_L(\tilde{X})$ ($\Pp^w_L(\tilde{X})$ resp.) the set of Pareto (weak Pareto resp.) optimal solutions in the sense of minimizing $L$ over $\tilde{X}$. W
    we have
    \begin{subequations}
    \begin{equation}\label{mono_x1}
    \Pp_L(\tilde{X}) \subseteq \Pp_L(X) \cup (\tilde{X} \setminus X), \quad \Pp^w_L(\tilde{X}) \subseteq \Pp^w_L(X) \cup (\tilde{X} \setminus X) \ .
    \end{equation}In particular,  
    \begin{equation}\label{mono_x2}
    \Pp_L(\tilde{X}) = \Pp_L(\Pp_L(X) \cup (\tilde{X} \setminus X)), \quad
    \Pp^w_L(\tilde{X}) = \Pp^w_L(\Pp^w_L(X) \cup (\tilde{X} \setminus X)) \ .
    \end{equation}
    \end{subequations}
    \end{enumerate}
\end{proposition}

\begin{remark}
In light of~\Cref{mono_proper}, the set of Pareto optimal meta-solvers is ``monotone'' with respect to additional performance criteria, that is, the candidates in the current Pareto optimal set will remain Pareto optimal as new criteria are introduced. 
On the other hand, when adding new meta-solvers, it is sufficient to compare the newly added solvers with the existing set of Pareto optimal meta-solvers.
\end{remark}
%The following property is interested for identifying the Pareto optimal solvers. 

\subsection{Re-scaling and Discovery of Optimal Solvers by Preference Functions}\label{subsec-discovery}
In this section, we propose a preference function based methodology for identifying the optimal solver among Pareto optimal ones. 
This approach encodes 
the features and properties that must be emphasized to select an optimal solver to solve a dedicated task in a specified context. 

First, the performance data for various criteria is evaluated using a re-scaling function. Note that other normalization methods can also be applied to the performance data, that map the original value in each dimension to a common subset of $\R$, allowing a unified evaluation. Here, without loss of generality, we use a re-scaling function. 
In particular, considering the MOO problem of minimizing~\eqref{MOO_obj} over $X$, for each dimension $i\in\{1,2,\dots,N\}$, denote
\begin{equation}
\overline{L}_i = \max_{x \in \Pp_L(X)} L_i(x), \quad \underline{L}_i = \min_{x \in \Pp_L(X)} L_i(x) \ ,
\end{equation}
and assume $\overline{L}_i \neq \underline{L}_i$ for every $i$. 
Then, for a solver $x_k \in \Pp_L(X)$, its performance is rescaled to
\begin{equation}
L'_i(x_k) = \frac{L_i(x_k) - \underline{L}_i}{\overline{L}_i - \underline{L}_i} \ .
\end{equation}

All the solvers on the set of Pareto optimal ones
cannot be dominated by any other solvers in terms of all criteria. However, in real application, users may have ``preference'' in favor of certain criteria over others. We model this by a preference function: $p: \R^N \to \R$, in the sense of minimizing the $N-$dimensional vector valued function $L: \R^d \to \R^N$. 
\begin{definition}\label{increasing_g}
For any $x^1,x^2 \in \R^N$, we denote $x^1 \leq x^2$ if $x^1_i \leq x^2_i$, for all $i \in \{1,2,\dots, N\}$. Moreover, we denote $x^1 < x^2$ if further there exists a $j$ such that $x^1_j < x^2_j$. 
A function $g: \R^N \to \R$ is called non-decreasing if for every $x^1, x^2 \in \R^N$ such that $x^1 \leq x^2$, we have $g(x^1) \leq g(x^2)$. A function $g: \R^N \to \R$ is called  increasing if for every $x^1, x^2 \in \R^N$ such that $x^1 < x^2$, we have $g(x^1) < g(x^2)$. 
\end{definition}
\begin{proposition}\label{propo_prefer}
Let $g: \R^N \to \R$ be an increasing function. Assume $\{L_i \}_{i \in \{ 1,2 \dots ,N\}}$ are lower-semicontinuous and proper functions. Then, 
\begin{equation} \text{if } \bar{x} \in \argmin_{x \in \R^d} g \circ L(x), \text{ we have } L(\bar{x}) \in \F_L(X) \ . \end{equation}
\end{proposition}

\begin{lemma}[Corollary of~\Cref{propo_prefer}]\label{lemma_prefer}
Let $p: \R^N \to \R$ be a  increasing function. Assume $\{L_i \}_{i \in \{ 1,2 \dots ,N\}}$ are lower-semicontinuous and proper functions. Then, we have 
\begin{equation} \argmin_{x \in X } g \circ L (x) = \argmin_{x \in \Pp_L(X)} g \circ L(x) \ . \end{equation}
\end{lemma}

As for identifying one optimal meta-solver considering the criteria $f$, that is the problem~\eqref{eq:moo}, 
given a preference function $p:\R^7 \to \R$ on all of the criteria, the solvers in the Pareto front will be evaluated again on the re-scaled data $f' = (f'_1,\dots, f'_7)$ by $p \circ (f'_1,f'_2,\dots,f'_7)$. That is, we solve a new classical optimization problem of the form:
\begin{equation}
\min_{x \in \Pp_L(X)} \ p \circ (f'_1,f'_2,\dots,f'_7) (x) \ .
\end{equation}
In this context, different users may employ their own preference functions, resulting in a plurality of optimal solvers. 

\subsection{Weighted-Sum Preference, Rediscovery of Optimal Solvers Based on Linear Programming}\label{sec-red}

\begin{figure}
    \subfigure[Sketch of the Pareto front of Meta-solvers on computational time and relative error.]{
    \includegraphics[width=0.45\textwidth]{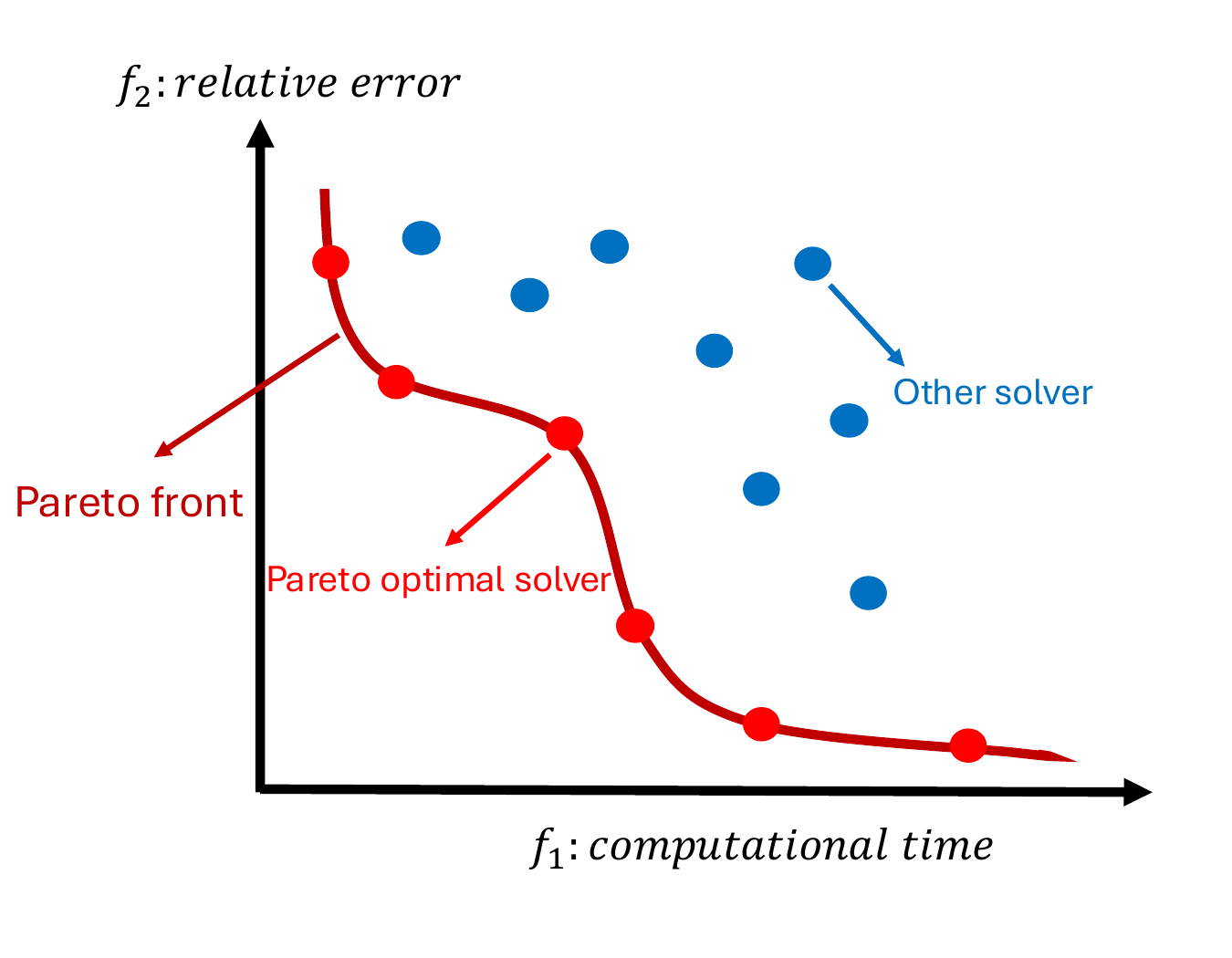}}
    \subfigure[Discovery of optimal meta-solvers by the weighed sum preference functions.]
    {\includegraphics[width=0.45\textwidth]{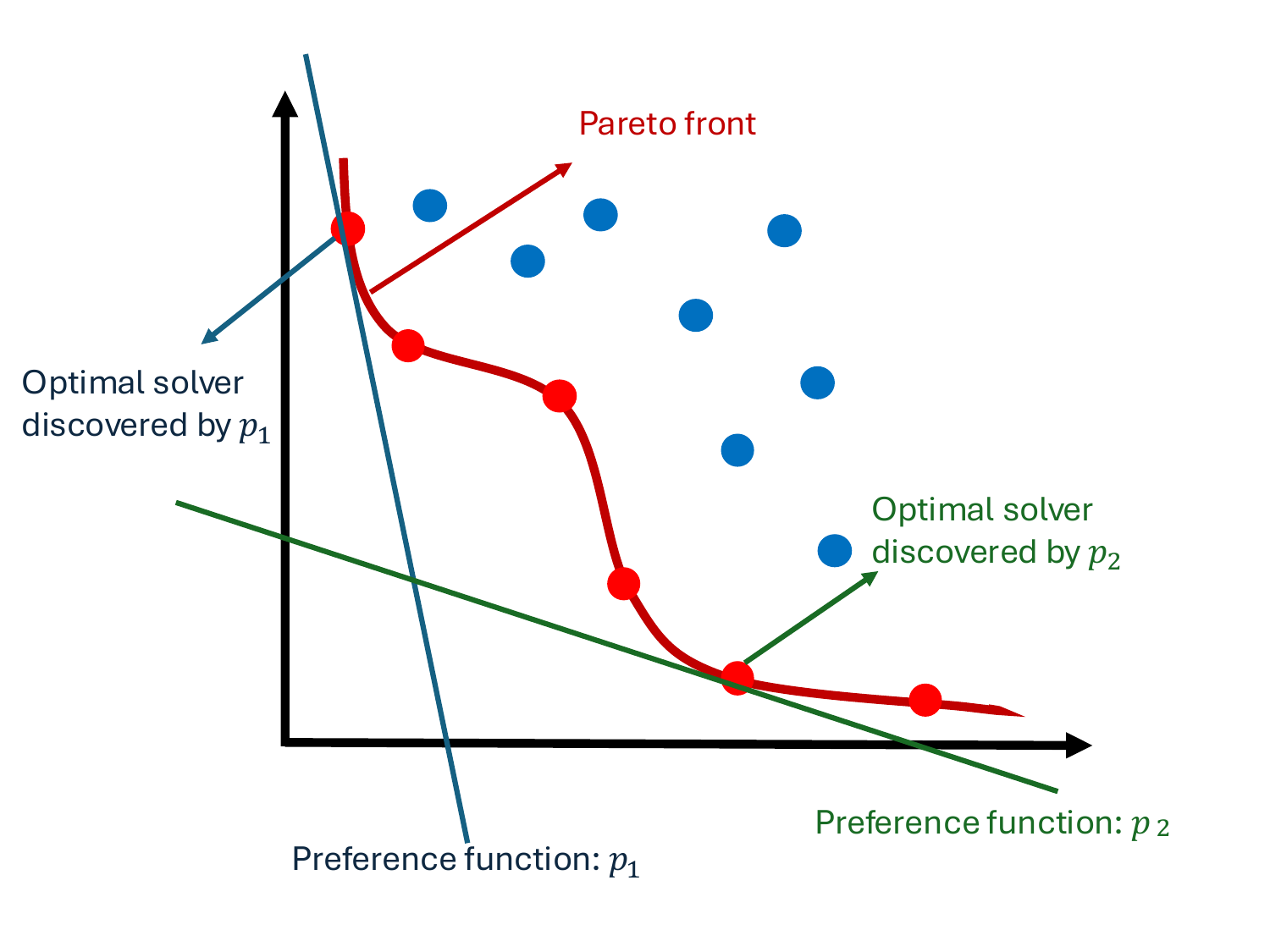}
    }
    \caption{Sketch of Pareto Front and discovering optimal solvers by preference functions, with computational time and relative error performance criteria. All solvers are depicted in blue while Pareto optimal solvers are highlighted in red. The discovery of optimal solvers by preference function geometrically corresponds to slicing the set of Pareto optimal meta-solvers with a straight line.  }
    \label{fig::discover}
\end{figure}

In this section, 
we explore the rediscovery of solvers based on a particular type of preference functions, the weighted sum function of the rescaled performance. Given a weight $\lambda = (\lambda_1,\lambda_2,\dots,\lambda_7)$ such that $0 \leq \lambda_i \leq 1, \forall \ i \in\{1,2,\dots, 7\}$ and $\sum_{i=1}^7\lambda_i = 1$, the preference function $p(\lambda; \cdot):\R^7 \to \R$ is defined by
\begin{equation}\label{weight_sum}
f'=(f'_1,f'_2,\dots,f'_7) \mapsto p(\lambda;f') := \sum_{i=1}^7 \lambda_i f'_i = \lambda^Tf' \ . 
\end{equation}
Given different weights $\lambda$, the preference function is used to select a solver among the Pareto optimal ones. 
Geometrically, each weighted sum preference function corresponds to a 7-dimensional hyperplane, where the weights determine the orientation of the hyperplane. One can incrementally decrease the constant term, starting from $0$, until the hyperplane intersects the Pareto front. The point of intersection corresponds to the performance of the optimal solver being selected, and the constant term represents the inverse of the weighted sum. 
\Cref{fig::discover} shows a sketch of the Pareto font of meta-solvers, and discovering optimal meta-solver by the weighted sum preference.  
%%%%%%%%%%%%%%%%%%%%%%

The re-discovery of one particular (parameterized) solver $x^j \in \Pp_{f'}(X)$ is equivalent to finding a weight $\lambda^j$ such that $p(\lambda;f'(x^j))$ is the minimum for every $x \in \Pp_{f'}(X)$. Here, $X \in \{\M^r,\M^K\}$ depends on the families of solvers. This can be formulated as a linear programming problem, that has the form:
\begin{equation}\label{lp_discover}
\begin{aligned}
&\min_{\lambda} \ \lambda^T f'(x^j) \\
&\text{ s.t. }
\left\{ 
\begin{aligned}
& \lambda^T f'(x^j) \leq \lambda^T f'(x^{-j}),\ \forall \ x^{-j} \in \Pp_X(f') \ ,\\
& 0 \leq \lambda_i \leq 1 \text{ and } \sum \lambda_i=1 \ .
\end{aligned}
\right.
\end{aligned}
\end{equation}
It is geometrically equivalent to finding the tangent hyperplane of the Pareto front at the point $f'(x^j)$. Note that the LP problem~\eqref{lp_discover} can be solved efficiently using modern LP solvers (e.g., CPLEX~\cite{cplex2009v12}).

%% file: numerical_poisson.tex
\section{Illustrative Numerical Examples: Solving Poisson equations}\label{sec-num}

In this section, as a concrete example for discovering the optimal meta-solver, we consider the Poisson equation with variable coefficient, which reads
\begin{equation}
\label{eqn:poisson}
\left\{
\begin{aligned}
    &- \nabla \cdot \left(\kappa \nabla u\right) = f, & \text{ in } \Omega \subset \R^{d} \ , \\
    &u = 0, & \text{ on } \partial\Omega \ ,
\end{aligned}
\right.
\end{equation}
where $\kappa$ and $f$ denote a coefficient function and a source function, respectively.
Note that the domain $\Omega$ is bounded in $\R^{d}$, with $d \in \{1, 2, 3\}$. 
We consider cases where equation~\eqref{eqn:poisson} has an unique solution. 
 %Implementation details can be found in \Cref{appd_imp}. 
 In the following, we present representative numerical results for each step of our methodology for both relaxation-based and Krylov-based methods. In the main paper, we present the results only in 3-d case, while the results in 1-d and 2-d cases can be found in~\Cref{app_numerics}.

\subsection{Numerical Results for Relaxation Based Meta-solvers}\label{subsec-relax}
In light of parameterization of meta-solvers in~\eqref{relax_para_solver}, we constructed and tested 588 different solvers for solving Poisson equation in 3-d, % in each of the three dimensions, 
i.e., $|\M^r| = 588$. 
%For solving $1-d$ Poisson equation, among 588 solvers we tested, there are 76 Pareto optimal solvers. For solving $2-d$ Poisson equation, t
There are 119 Pareto optimal solvers.

\subsubsection{Identifying Optimal Solvers in Pareto Sense}
We present the composition of the Pareto optimal solvers, by counting the number of different $x_i$. The results are shown in~\Cref{relax_com}. As a preliminary finding, we observe that JacobiKAN accounts for the largest number of Pareto optimal solvers. Among classical iterative solvers, SSOR produces the largest and most dominant count. The 2-level multi-grid technique accounts for the majority, followed by the 3-level multi-grid technique. Finally, the most efficient ratio is $\frac{1}{128}$, the smallest one, and the number decreases as the ratio increases. This trend is consistent with the fact that in 3-d problems, there are more high-frequency components that need to be relaxed, which requires a higher number of iterations for the relaxation methods. 
%In particular, The number of different neural operators is shown in~\Cref{n_no_pareto}, the number of different classical iterative solvers is shown in~\Cref{n_co_pareto}, the number of different strategies (different levels) for adapting multi-grids techniques is shown in~\Cref{n_mt_pareto}, the number of different ratio of neural operators is shown in~\Cref{n_ra_pareto}.

\begin{table}
    \centering
    \small
    \caption{The composition of the set of Pareto optimal solvers by counting the number of elements in each dimension for relaxation based methods, solving 3-d Poisson equation, with the highest number highlighted.}
    \begin{minipage}{\linewidth}
        \centering
        \captionof{subtable}{Different neural operators,}
        \resizebox{\linewidth}{!}{ 
        \begin{tabular}{c| c c c c c c c} 
\hline
Neural Op & DeepONet & U-Net & FNO & Transformer & KAN & JacobiKAN& ChebyKAN  \\
\hline
%\# in Pareto opt in 1D  & 7 & 5 & 11 & 11 & 13 & \textbf{17} & 12 \\
%\hline
%\# in Pareto opt in 2D  & 16 & \textbf{18} & 4 & 5 & 15 & 17 & 0 \\
%\hline
\# in Pareto opt  & 20 & 18 & 7 & 5 & 21 & \textbf{27} & 21 \\
\hline
%\# in Pareto opt., $\frac{1}{128}$ & 28 & 20 & 7 & 8 & 14 & 36
\end{tabular}
}
    \end{minipage}

    \vspace{1em} % Space between subtables

    % Second subtable
    \begin{minipage}{\linewidth}
        \centering
        \captionof{subtable}{Different classical iterative solvers.}
           \begin{tabular}{c| c c c c    } 
\hline
Classical solvers & Jacobi & Gauss-Siedel & SOR & SSOR  \\
\hline
%\# in Pareto opt in 1D & 2 & 19 & 22 & \textbf{33}  \\
%\hline
%\# in Pareto opt in 2D & 2 & 12 & 23 & \textbf{38}  \\
%\hline
\# in Pareto opt & 8 & 14 & 37 & \textbf{60}  \\
\hline
%\# in Pareto opt., $\frac{1}{128}$ & 1 & 8 & 26 & 27 & 35 & 16 & 0 
\end{tabular}
        
    \end{minipage}

\begin{minipage}{\linewidth}
        \centering
        \captionof{subtable}{Different number of multi-grids.}
           \begin{tabular}{c| c c c  } 
\hline
Multi-grid method& NO & 2-level & 3-level \\
\hline
%\# in Pareto opt in 1D & 1 & \textbf{56} & 19 \\
%\hline
%\# in Pareto opt in 2D & 13 & \textbf{43} & 19 \\
%\hline
\# in Pareto opt & 15 & \textbf{63} & 41 \\
\hline
%\# in Pareto opt., $\frac{1}{128}$ & 62 & 42 & 9 \\
\end{tabular}
        
    \end{minipage}

\begin{minipage}{\linewidth}
        \centering
        \captionof{subtable}{Different proportion of neural operators.}
    \begin{tabular}{c| c c c c c c c } 
\hline
Ratio of N.O. & $\frac{1}{2}$ & $\frac{1}{4}$ & $\frac{1}{8}$ & $\frac{1}{16}$ & $\frac{1}{32}$ & $\frac{1}{64}$ & $\frac{1}{128}$  \\ 
\hline
%\# in Pareto opt in 1D & 0 & 0 & 5 & 13 & {13} & {21} & \textbf{24} \\ 
%\hline
%\# in Pareto opt in 2D & 0 & 0 & 3 & 9 & 14 & 24 & \textbf{25} \\ 
%\hline
\# in Pareto opt & 0 & 2 & 12 & 18 & 27 & 28 & \textbf{32} \\ 
\hline
%\# in Pareto opt., $\frac{1}{128}$ & 6 & 6 & 10 & 22 & 23 & 23 & 23 \\ 
\end{tabular}

       \end{minipage}
    \label{relax_com}
\end{table}

\begin{figure}
\centering
\subfigure[Computational time -- Relative error -- \# of iterations]{
\includegraphics[width=0.3\textwidth]{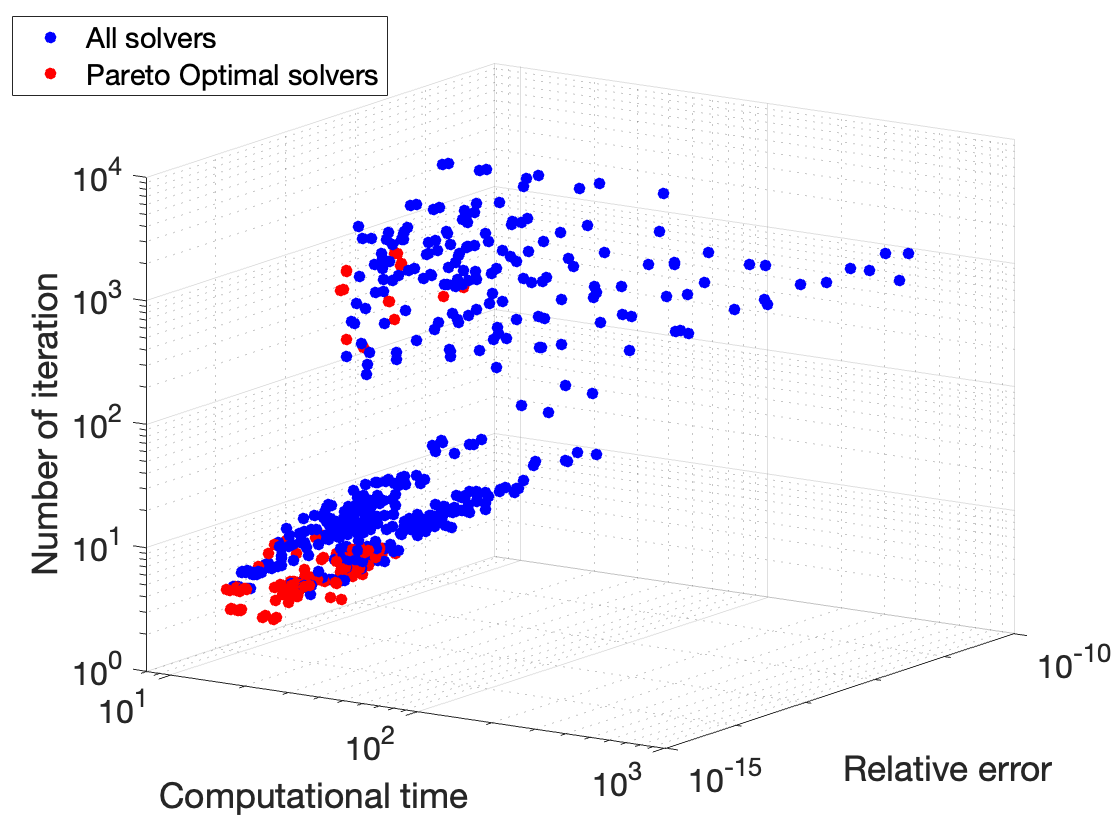}}
\subfigure[Relaxtive error -- MACs -- Memory allocation]{
\includegraphics[width=0.3\textwidth]{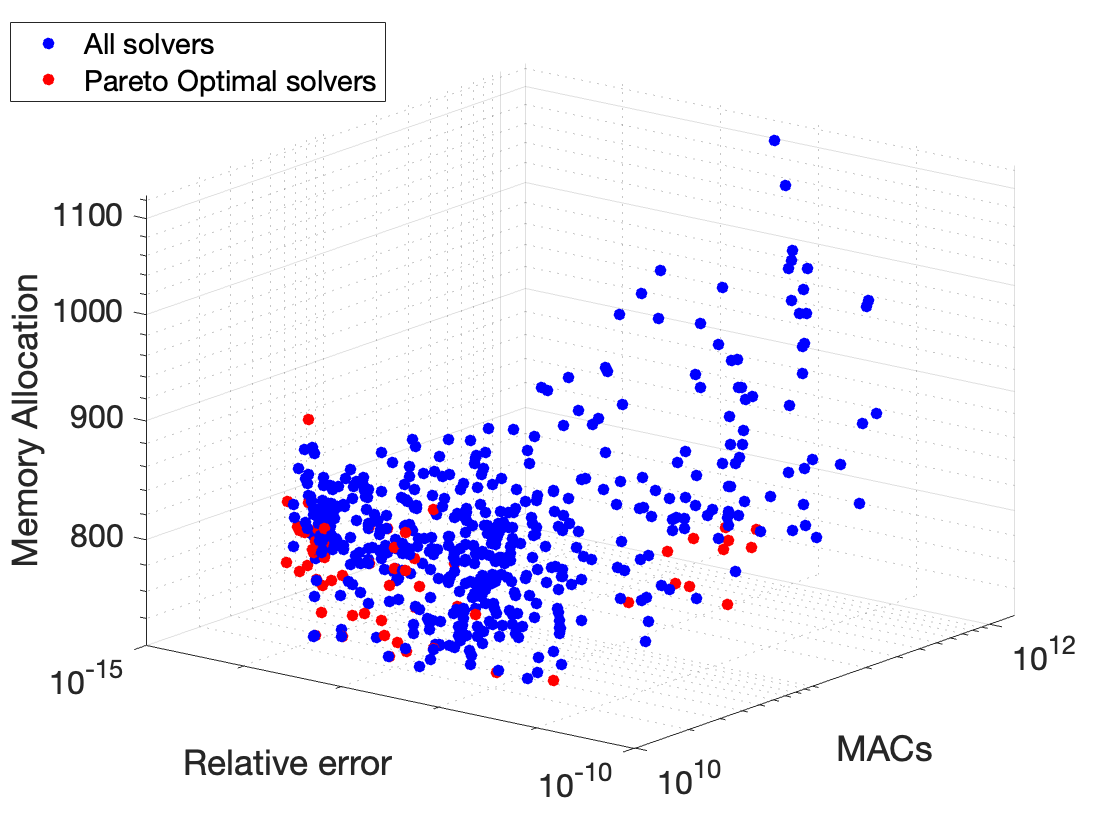}}
\subfigure[Convergence rate -- MACs -- \# of iterations]{
\includegraphics[width=0.3\textwidth]{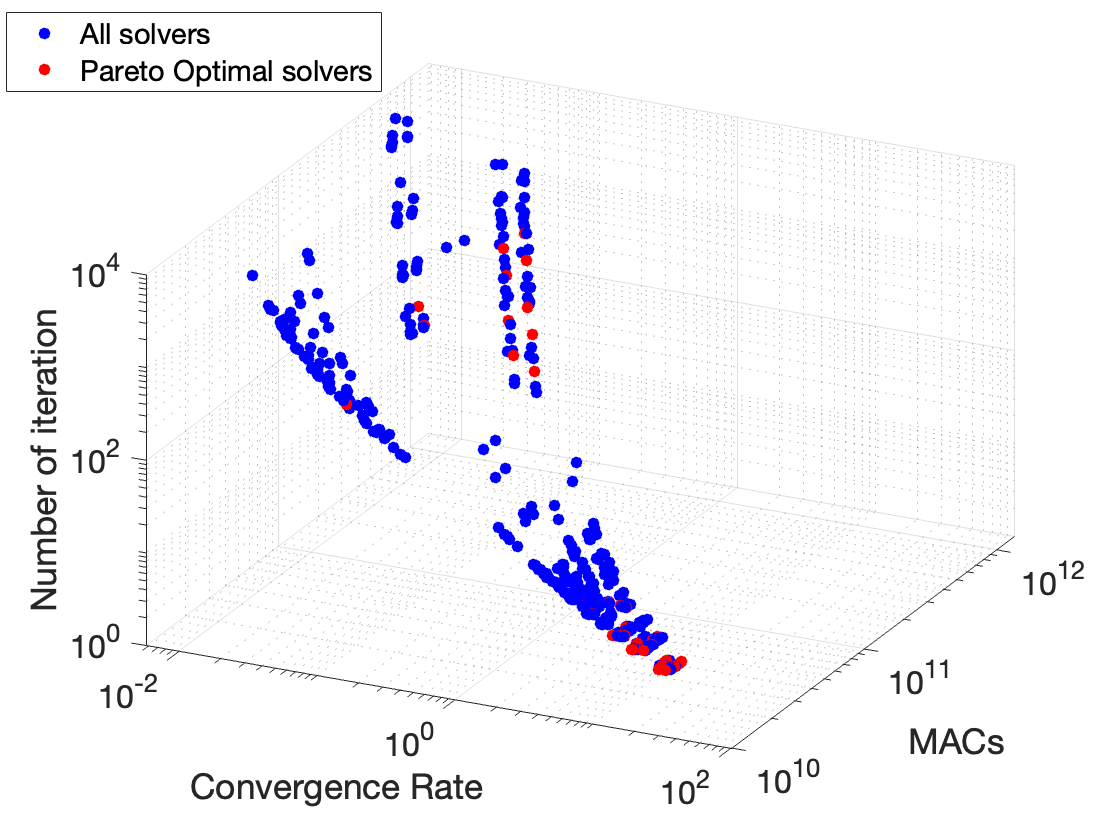}}
\caption{Projection of Pareto front intro three dimensional criteria for solving 3-d Poisson equation, using ralxation-based methods. All solvers are depicted in blue while Pareto optimal solvers are highlighted in red. The ``gap'' due to the adaption thus the improvement of performance by multi-grid techniques.}
\label{front_3d_relax}
\end{figure}

We also plot the projection of Pareto fronts onto various 3-dimensional criteria in~\Cref{front_3d_relax}. These figures show that most performance criteria provide significant variability that allows us to discriminate between meta-solvers. However, we also observe that some performance criteria are highly correlated, e.g., convergence rate and number of iterations. Note that this observation does not prevent the discovery of optimal meta-solvers. 
 These results also depict visualization of the projection of the Pareto front in various three-dimensional representations to obtain numerical evidence on the topology and shape of the Pareto front. 

%for solving $1-d, 2-d$ and $3-d$ Poisson equations. 
%First, we show the projections for computational time, relative error, and number of iterations in~\Cref{front_tei_12}. Next, we present the projections for relative error, MACs, and memory allocation in~\Cref{front_emm_12}, and finally, for convergence rate, MACs, and number of iterations in~\Cref{front_rmi_12}.

%\begin{figure}
%\centering
%\subfigure[Solving 1d Poisson equation]{
%\includegraphics[width=0.3\textwidth]{figs/1d_emm.png}}
%\subfigure[Solving 2d Poisson equation]{
%\includegraphics[width=0.3\textwidth]{figs/2d_emm.png}}
%\subfigure[Solving 3d Poisson equation]{
%\includegraphics[width=0.3\textwidth]{figs/3d_emm.png}}
%\caption{Projection of Pareto front intro three dimensional criteria: Relative Error -- MACs -- Memory Allocation, for relaxation-based methods.}
%\label{front_emm_12}
%\end{figure}

%\begin{figure}
%\centering
%\subfigure[Solving 1d Poisson equation]{
%\includegraphics[width=0.3\textwidth]{figs/1d_rmi.png}}
%\subfigure[Solving 2d Poisson equation]{
%\includegraphics[width=0.3\textwidth]{figs/2d_rmi.png}}
%\subfigure[Solving 3d Poisson equation]{
%\includegraphics[width=0.3\textwidth]{figs/3d_rmi.png}}
%\caption{Projection of Pareto front intro three dimensional criteria: Convergence Rate -- MACs -- \# of iterations, for relaxation-based methods.}
%\label{front_rmi_12}
%\end{figure}

\subsubsection{Discovery of Optimal Solvers via Preference}\label{subsec_dis}
Different users may employ their own preference functions, resulting in various optimal solvers. To highlight the realtive performance, we present the top three solvers for specific preference functions applying on the re-scaling data. %, in solving 1-D, 2-D, and 3-D Poisson equations.

\begin{pref}
$p^1(f') = \frac{1}{7}(\sum_{i=1}^7 f'_i)$. The average of all the re-scaling performance. 
\end{pref} 
%,  solving $2-d$ Poisson equation is shown in~\Cref{top3_p1_2d}, 
%and solving $3-d$ Poisson equation is shown in~\Cref{top3_p1_3d}. 
\begin{pref}
 $p^2(r) =  0.4 \times (r_1 + r_5) + 0.04 \times (r_2 + r_3 + r_4 + r_6 + r_7)$. This preference function means that the computational time and relative error are ten times more important then other criteria. 
 \end{pref}
 
\begin{figure}
\centering
\subfigure[Top-3 solvers discovered by preference function $p^1$]{
\includegraphics[width=0.45\textwidth]{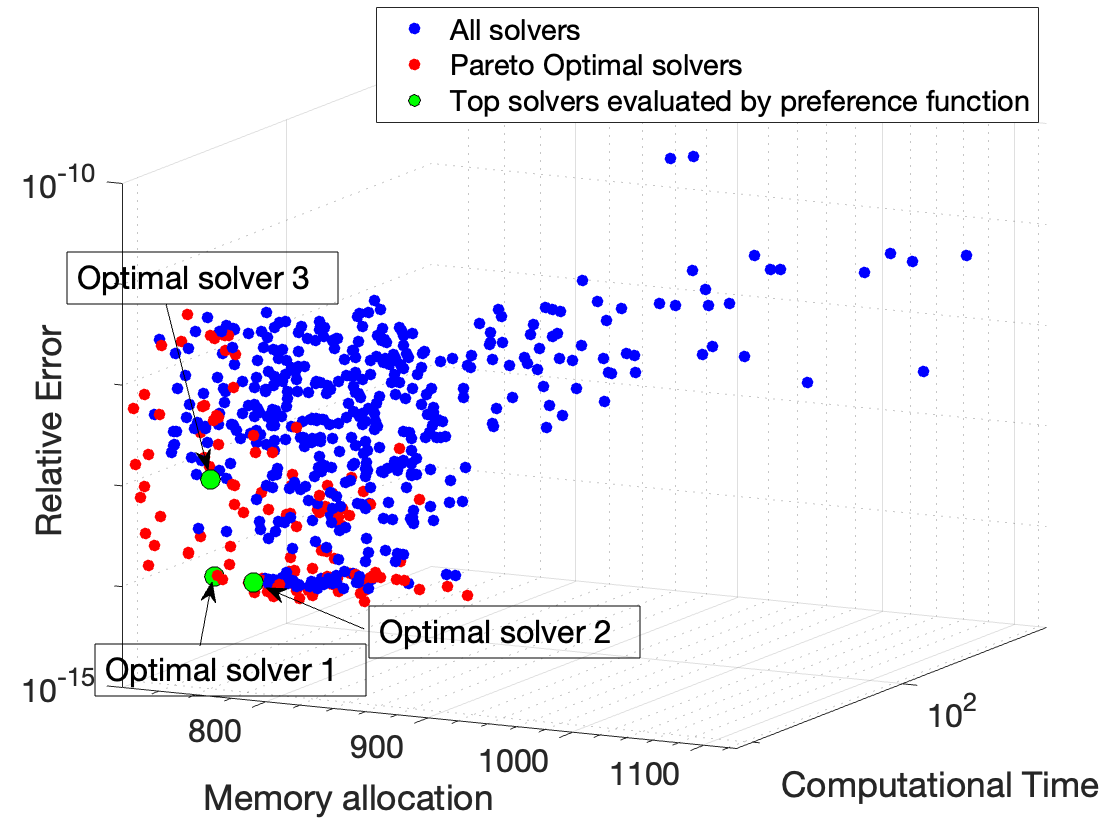}}
\subfigure[Top-3 solvers discovered by preference function $p^2$]{
\includegraphics[width=0.45\textwidth]{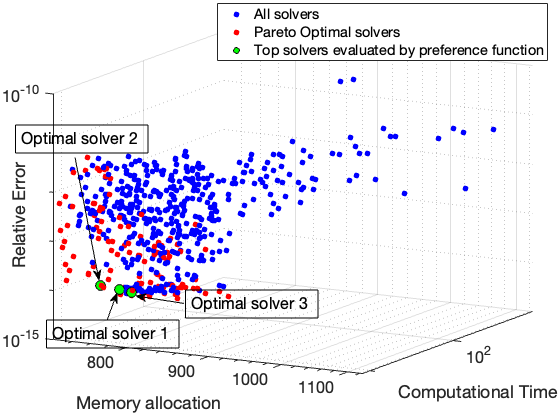}}
\caption{The top-3 optimal solvers discovered by preference function $p^1$ and $p^2$, in the three-dimensional projection: Computational Time -- Relative Error -- Memory allocation of Pareto front, for relaxation-based methods, solving 3-d Poisson equation.}
\label{front_dis_relax}
\end{figure}

The top 3 solvers and their performance using preference function $p^1$ for solving the 3-d Poisson equation is shown in~\Cref{top3_p1_3d}, while the results % The top 3 solvers and their performance 
using preference function $p^2$ is shown in~\Cref{top3_p2_3d}. 
We observe that some consistency is maintained among the optimal solvers identified by the two distinct preference functions. In particular, the top two solvers are the same, though their order differs. The best ratio is always $\frac{1}{128}$. DeepONet emerges as the dominant neural operator (appearing in 5 out of 6 cases), while SSOR is the dominant classical iterative solver (also appearing in 5 out of 6 cases). These combinations, along with the resulting meta-solvers, are optimal both in terms of average performance across the criteria and in terms of speed and accuracy. Moreover, we plot in \Cref{front_dis_relax} the top-3 solvers discovered by both preference functions on the Pareto front, considering the criteria of computational time, relative error, and memory allocation. The positions in the performance figure of the discovered optimal solvers align with the observation found in the tables.

\begin{table}
    \centering
    \caption{Top-3 solvers by preference function $p^1$ in 3-D for relaxation-based methods.}
    
    \begin{minipage}{\linewidth}
        \centering
        \captionof{subtable}{The top-3 solvers.}
        \resizebox{\linewidth}{!}{
        \begin{tabular}{c|c|c|c|c}
\hline \hline 
 &Neural operator & Classical solver & Multi-grid &  Ratio \\ 
\hline \hline
Top 1 solver of $p^1$ (3d t1\_1)  & DeepONet & SSOR & 3-level & $\frac{1}{128}$ \\ \hline
Top 2 solver of $p^1$ (3d t2\_1)  & DeepONet & SSOR & 2-level &  $\frac{1}{128}$ \\ \hline
Top 3 solver of $p^1$ (3d t3\_1)  & DeepONet & SOR & 3-level &  $\frac{1}{128}$ \\ \hline
%Accurate-4($\frac{1}{128}$)& ChebyKAN & SSOR & 1-level & 0 & 0.015625 \\ \hline
%Accurate-5($\frac{1}{128}$)& DeepONet & BiCGStab & 1-level & 8 & 0.5 \\  \hline
\end{tabular}
}
    \end{minipage}

    \vspace{1em} % Space between subtables

    % Second subtable
    \begin{minipage}{\linewidth}
        \centering
        \captionof{subtable}{Performance of the top-3 solvers.}
        \resizebox{\linewidth}{!}{ % Resize to fit the line width
          \begin{tabular}{c|c|c|c|c|c|c|c}
\hline \hline 
 &Error & Com. time &  \# of ite. & Conv. rate & Memory & MACs & Training time \\ 
\hline \hline
3d t1\_1& 1.42E-14 & 9.047 & 2 & 17.798 & 765.639 & 2.291E+10 & 8625.91 \\ \hline
3d t2\_1& 1.41E-14 & 8.308 & 2 & 17.851 & 790.796 & 2.274E+10 & 8625.91 \\ \hline
3d t3\_1& 1.31E-13 & 8.880 & 2 & 14.792 & 764.098 & 2.244E+10 & 8625.91 \\ \hline
%Accurate-4($\frac{1}{128}$)& 0.0971029 & 1.24E-14 & 6 & 2.86298 & 235.11 & 178858 & 99.2164 \\ \hline
%Accurate-5($\frac{1}{128}$)& 0.0518639 & 1.32E-14 & 13 & 1.47149 & 236.639 & 178858 & 326.378 \\  \hline
\end{tabular}
        }
    \end{minipage}
    \label{top3_p1_3d}
\end{table}

\begin{table}
\small
    \centering
    \caption{Top-3 solvers by preference function $p^2$ in 3-D for relaxation-based method.}
    \begin{minipage}{\linewidth}
        \centering
        \captionof{subtable}{The top-3 solvers.}
        \resizebox{\linewidth}{!}{ 
        \begin{tabular}{c|c|c|c|c}
\hline \hline 
 &Neural operator & Classical solver & Multi-grid &  Ratio \\ 
\hline \hline
Top 1 solver of $p^2$ (t1\_2 3d) & DeepONet & SSOR & 2-level & $\frac{1}{128}$ \\ \hline
Top 2 solver of $p^2$ (t2\_2 3d) & DeepONet & SSOR & 3-level &  $\frac{1}{128}$ \\ \hline
Top 3 solver of $p^2$ (t3\_2 3d) & FNO & SSOR & 2-level &  $\frac{1}{128}$ \\ \hline
%Accurate-4($\frac{1}{128}$)& ChebyKAN & SSOR & 1-level & 0 & 0.015625 \\ \hline
%Accurate-5($\frac{1}{128}$)& DeepONet & BiCGStab & 1-level & 8 & 0.5 \\  \hline
\end{tabular}
}
    \end{minipage}

    \vspace{1em} % Space between subtables

    % Second subtable
    \begin{minipage}{\linewidth}
        \centering
        \captionof{subtable}{Performance and rank of performance of the top-3 solvers.}
        \resizebox{\linewidth}{!}{ % Resize to fit the line width
           \begin{tabular}{c|c|c|c|c|c|c|c}
\hline \hline 
 &Error & Com. time &  \# of ite. & Conv. rate & Memory & MACs & Training time \\ 
\hline \hline
t1\_2 3d & 1.41E-14 & 8.308 & 2 & 17.851 & 790.796 & 2.274E+10 & 8625.9 \\ \hline
t2\_2 3d & 1.42E-14 & 9.047 & 2 & 17.798 & 765.639 & 2.291E+10 & 8625.9 \\ \hline
t3\_2 3d & 1.34E-14 & 8.354 & 2 & 17.151 & 804.265 & 2.137E+10 & 47736.7 \\ \hline
%Accurate-4($\frac{1}{128}$)& 0.0971029 & 1.24E-14 & 6 & 2.86298 & 235.11 & 178858 & 99.2164 \\ \hline
%Accurate-5($\frac{1}{128}$)& 0.0518639 & 1.32E-14 & 13 & 1.47149 & 236.639 & 178858 & 326.378 \\  \hline
\end{tabular}
        }
    \end{minipage}
    \label{top3_p2_3d}
\end{table}

\subsubsection{Rediscovery of Optimal Solvers Based on Linear Programming}

We implement the method introduced in~\Cref{sec-red} for re-discovering particular solvers. 
In particular, given a solver and its performance, the weights in the preference function are computed by solving the problem~\eqref{lp_discover}. 
In the following, we present the weights in the order s.t. $(\lambda_1,\lambda_2,\dots,\lambda_7)$ for ( $\lambda_1$: relative error, $\lambda_2$: computational time, $\lambda_3$: number of iterations, $\lambda_4$: convergence rate, $\lambda_5$: memory allocation, $\lambda_6$: MACs, $\lambda_7$: training time), for some particular solvers. 

\begin{subequations}\label{redis_lambda}
\begin{equation}\label{redis_a}
\left\{
\begin{aligned}
&\text{Preference function $p(\lambda_a;r) = (0.087, 0, 0, 0, 0.004, 0.861,  0.046)^Tr$ } \\
&\text{Optimal solver}: x^a = (\text{ChebyKAN}, \text{SSOR}, \text{2-level}, \frac{1}{128}) \ .
\end{aligned}
\right.
\end{equation}
\begin{equation}
\left\{
\begin{aligned}
&\text{Preference function $p(\lambda_b;r) = (0, 0, 0, 0.266, 0.734, 0,  0)^Tr$ } \\
&\text{Optimal solver}: x^b = (\text{DeepONet}, \text{SSOR}, \text{3-level}, \frac{1}{128}) \ .
\end{aligned}
\right.
\end{equation}
%\begin{equation}
%\left\{
%\begin{aligned}
%&\text{Preference function $p(\lambda_c;r) = (0.102, 0.039, 0, 0.067, 0.202, 0.590,  0)^Tr$ } \\
%&\text{Optimal solver}: x^c = (\text{Transformer}, \text{SOR}, \text{2-level}, \frac{1}{64}) \ .
%\end{aligned}
%\right.
%\end{equation}
\end{subequations}

As seen in \eqref{redis_lambda}, when the optimal solver changes, the weights in the preference function used for discovery also change significantly. 
This is consistent with the observation that the meta-solver within the Pareto optimal candidate is only optimal in certain contexts.
For instance, consider solver $x^a$ re-discovered in~\eqref{redis_a},  where the user places the most importance on memory allocation (with a weight of 0.861) and relative error (weight of 0.087), while assigning minimal importance to the other criteria (weights $<$ 0.05). 
 In this case, the meta-solver ($x_1$=ChebyKAN, $x_2$=SSOR, $x_3$=2-level, $x_4=\frac{1}{128}$) is the optimal one.

\subsection{Numerical Results for Krylov Based Meta-solvers}\label{sec-poi-kry}
In light of the parameterization of meta-solvers in~\eqref{krylov_para_solver}, we constructed 900 different Krylov based meta-solvers for the Poisson equation in 3-d, i.e., $|\M^K| = 900$. % for $d=1,2,3$. 
Among them there are 210 Pareto optimal solvers.

\subsubsection{Identifying Optimal Solvers  in the Pareto Sense}

We present the composition of the Pareto optimal solvers, by counting the number of different $x_i$. The results are shown in~\Cref{kry_com}. 
%In particular, The number of different neural operators is shown in~\Cref{n_no_pareto}, the number of different classical iterative solvers is shown in~\Cref{n_co_pareto}, the number of different strategies (different levels) for adapting multi-grids techniques is shown in~\Cref{n_mt_pareto}, the number of different ratio of neural operators is shown in~\Cref{n_ra_pareto}.
As a preliminary finding, we observe that DeepONet is the most dominant neural operator. Among the Krylov methods, FBiCGStab is the most dominant, followed by FGMRES. SSOR is the best smoother, accounting for over 65\% of all smoothers, while the proportions of different strategies applying smoothers are nearly identical. The 2-level multi-grid technique is the most prevalent, followed by the 3-level technique.

%Among classical iterative solvers, SSOR produces the largest and most dominant count. The 2-level multi-grid technique accounts for the majority, followed by the 3-level technique. Finally, the most efficient ratio is $\frac{1}{128}$, the smallest one, and the number decreases as the ratio increases. This trend is consistent with the fact that in 3D problems, there are more high-frequency components that need to be relaxed, which requires a higher number of iterations for the relaxation methods. 

\begin{table}
    \centering
    \small
    \caption{The composition of the set of Pareto optimal solvers by counting the number of elements in each dimension for Krylov based methods, solving 3-d Poisson equation, with the highest number highlighted.}
    \begin{minipage}{\linewidth}
        \centering
        \captionof{subtable}{Different neural operators.}
        \resizebox{\linewidth}{!}{ 
        \begin{tabular}{c| c c c  c c } 
\hline
Neural Op & DeepONet & U-Net  & KAN & JacobiKAN& ChebyKAN \\
\hline
%\# in Pareto opt in 1-d & 27 & 5 & 11 & \textbf{29} & 8 \\
%\hline
%\# in Pareto opt in 2-d & 25 & 42 & 56 & \textbf{62} & 36 \\
%\hline
\# in Pareto opt in 3-d & \textbf{61} & 31 & 39 & 28 & 51 \\
\hline
%\# in Pareto opt., $\frac{1}{128}$ & 28 & 20 & 7 & 8 & 14 & 36
\end{tabular}
}
    \end{minipage}

    \vspace{1em} % Space between subtables

    % Second subtable
    \begin{minipage}{\linewidth}
        \centering
        \captionof{subtable}{Different Krylov solvers.}
           \begin{tabular}{c| c c c c    } 
\hline
Classical solvers & FGMRES & FCG & FBiCGStab \\
\hline
%\# in Pareto opt in 1-d & 22 & 13 & \textbf{45}  \\
%\hline
%\# in Pareto opt in 2-d & 53 & 50 & \textbf{118}  \\
%\hline
\# in Pareto opt in 3-d & 61 & 39 & \textbf{110}  \\
\hline
%\# in Pareto opt., $\frac{1}{128}$ & 1 & 8 & 26 & 27 & 35 & 16 & 0 
\end{tabular}
        
    \end{minipage}

    \begin{minipage}{\linewidth}
        \centering
        \captionof{subtable}{Different smoothers.}
           \begin{tabular}{c| c c c c    } 
\hline
Smoother & GS & Jacobi & SOR & SSOR  \\
\hline
%\# in Pareto opt in 1-d & 15 & 15 & 15 & \textbf{35}  \\
%\hline
%\# in Pareto opt in 2-d & 43 & 39 & 56 & \textbf{83}  \\
%\hline
\# in Pareto opt in 3-d & 16 & 7 & 49 & \textbf{138}  \\
\hline
%\# in Pareto opt., $\frac{1}{128}$ & 1 & 8 & 26 & 27 & 35 & 16 & 0 
\end{tabular}
        
    \end{minipage}

\begin{minipage}{\linewidth}
        \centering
        \captionof{subtable}{Different strategies of applying smoothers.}
          \begin{tabular}{c| c c c c c    } 
\hline
Classical solvers & 1-1-1 & 3-1-3 & 5-1-5 & 7-1-7 & 9-1-9  \\
\hline
%\# in Pareto opt in 1-d & 20 & 6 & 10 & 19 & \textbf{25} \\
%\hline
%\# in Pareto opt in 2-d & 16 & \textbf{57} & 41 & 50 & \textbf{57} \\
%\hline
\# in Pareto opt in 3-d & 34 & 43 & 44 & \textbf{54} & 35 \\
\hline
%\# in Pareto opt., $\frac{1}{128}$ & 1 & 8 & 26 & 27 & 35 & 16 & 0 
\end{tabular}
        
    \end{minipage}

\begin{minipage}{\linewidth}
        \centering
        \captionof{subtable}{Different levels of multigrid.}
    \begin{tabular}{c| c c c     } 
\hline
Levels in multi-grid & 1-level & 2-level &  3-level   \\
%\hline
%\# in Pareto opt in 1-d & \textbf{46} & 28 & 6  \\
%\hline
%\# in Pareto opt in 2-d & 68 & \textbf{88} & 65  \\
\hline
\# in Pareto opt in 3-d & 57 & \textbf{83} & 70  \\
\hline
%\# in Pareto opt., $\frac{1}{128}$ & 1 & 8 & 26 & 27 & 35 & 16 & 0 
\end{tabular}

       \end{minipage}
    \label{kry_com}
\end{table}

\begin{figure}
\centering
\subfigure[Computational time -- Relative error -- \# of iterations]{
\includegraphics[width=0.3\textwidth]{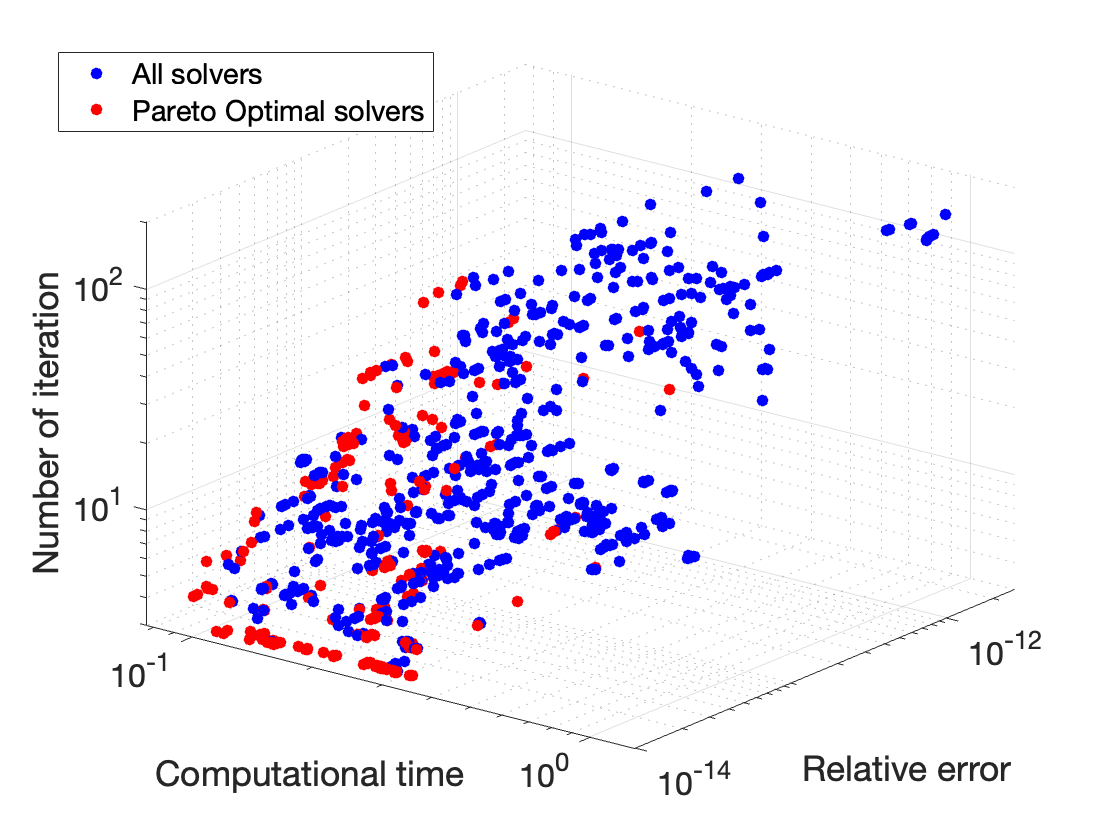}}
\subfigure[Relative error -- MACs -- Memory allocation]{
\includegraphics[width=0.3\textwidth]{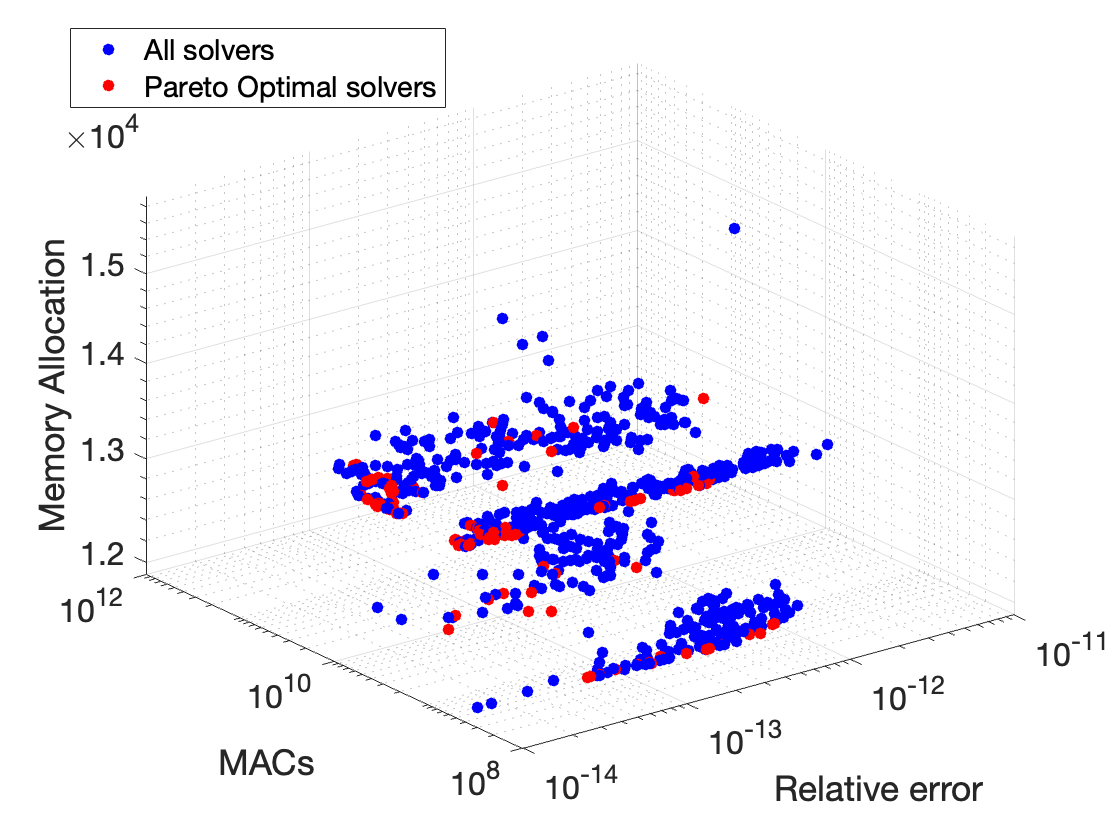}}
\subfigure[Convergence rate -- MACs -- \# of iterations]{
\includegraphics[width=0.3\textwidth]{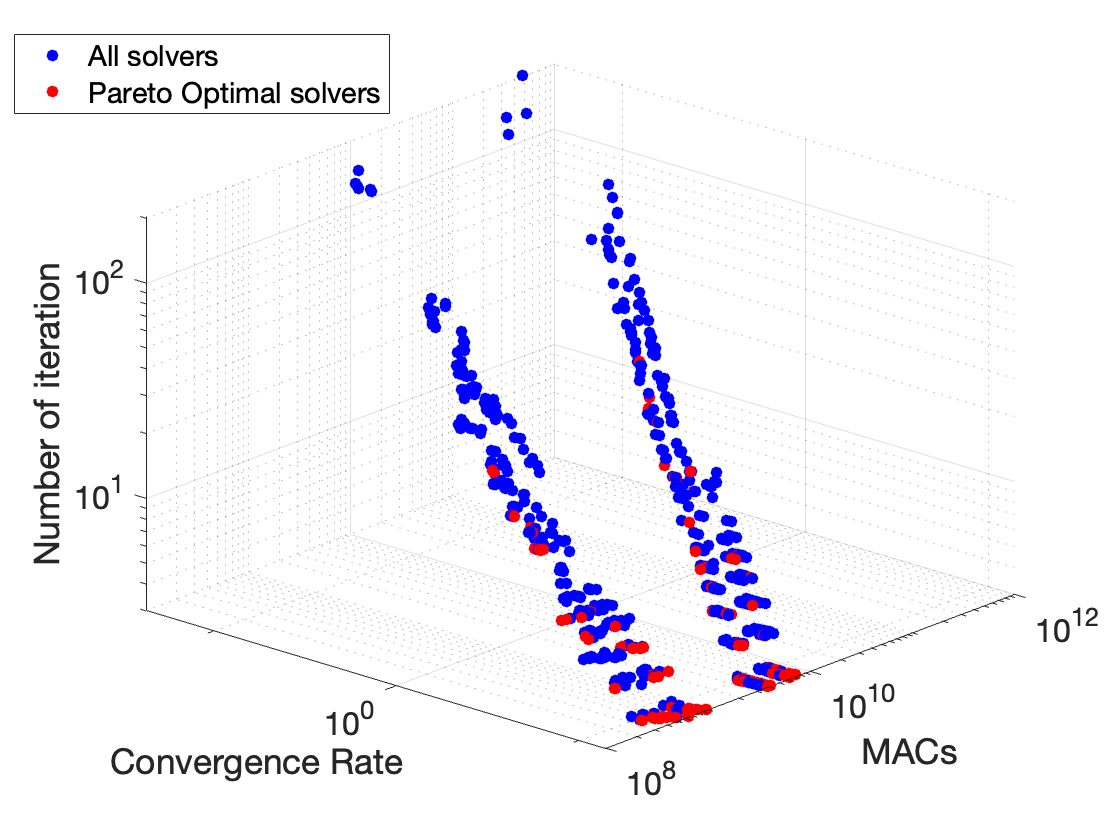}}
\caption{Projection of Pareto front onto three dimensional criteria for solving 3-d Poisson equation, using Krylov methods. All solvers are depicted in blue while Pareto optimal solvers are highlighted in red. The``gap'' due to the adaption thus the improvement of performance by multi-grid techniques.}
\label{front_3d_kry}
\end{figure}

We also plot the projection of the Pareto fronts onto various 3-dimensional criteria in \Cref{front_3d_kry}. 
Compared to the Pareto fronts of the relaxation-based methods shown in \Cref{front_3d_relax}, although the proportion of Pareto optimal solvers among all solvers is similar, we observe that the performance of different Krylov-based meta-solvers exhibits greater variability. 
This suggests a greater demand for methodologies in the discovery and re-discovery of optimal solvers.
Additionally, we observe that the shape and topology of the Pareto front for Krylov-based methods are significantly more complex than those for relaxation-based methods.

\subsubsection{Discovery/Re-discovery of Optimal Solvers}
\begin{figure}
\centering
\subfigure[Solving 1d Poisson Equation]{
\includegraphics[width=0.45\textwidth]{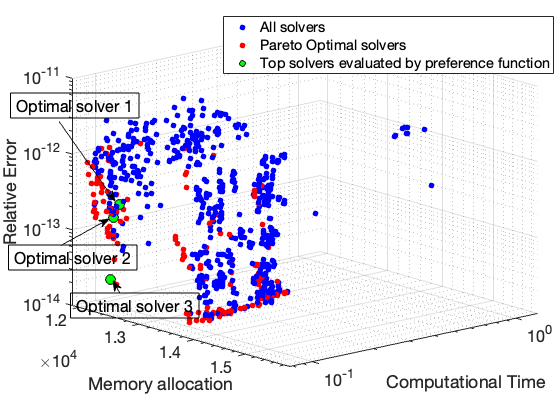}}
\subfigure[Solving 2d Poisson Equation]{
\includegraphics[width=0.45\textwidth]{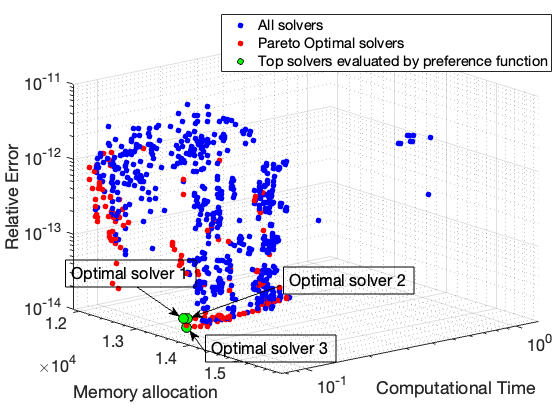}}
\caption{The top-3 optimal solvers discovered by preference functions $p^1$ and $p^2$, in the three-dimensional projection: Computational Time -- Relative Error -- Memory allocation of Pareto front, for Krylov-based methods, solving 3-d Poisson equation.}
\label{front_kry}
\end{figure}

We apply the same preference function $p^1$ and $p^2$ as for relaxation based methods in~\Cref{subsec_dis}. 
For $p^1$, the top-3 optimal solvers discovered is shown in~\Cref{top3_p1_3d_kry}. For $p^2$, the top-3 optimal solvers discovered is shown in~\Cref{top3_p2_3d_kry}.  
We plot the top-3 solvers in the Pareto front, for criteria computational time, relative error and memory allocation, in~\Cref{front_kry}. % in 1-D, 2-D and 3-D cases. 
We observe that different preference functions lead to much different optimal solvers, although some consistency is maintained. This is consistent with the observation that Krylov-based meta-solvers exhibit greater variability compared to relaxation-based meta-solvers. 
In particular, FBiCGStab is consistently the best Krylov method, and SSOR is consistently the best relaxation method, under both preference functions $p^1$ and $p^2$. However, when considering the average performance across all criteria, the 1-level multi-grid method is the best, while the 3-level multi-grid method excels when the focus is on computational time and relative error. 
This distinction is clearly illustrated in~\Cref{front_kry}, where the discovered optimal solvers by two preference functions are in two different ``point sets'' of solvers, corresponding to performance of using 1-level multi-grid and 3-level multi-grid.

\begin{table}
    \centering
    \caption{Top-3 solvers by preference function $p^1$ in 3-D for Krylov-based methods.}
    
    \begin{minipage}{\linewidth}
        \centering
        \captionof{subtable}{The top-3 solvers.}
        \resizebox{\linewidth}{!}{
        \begin{tabular}{c|c|c|c|c|c}
\hline \hline 
 &Neural operator & Classical solver & Multi-grid & Relaxation & Strategies \\ 
\hline \hline
Top 1 solver of $p^1$ (3d t1\_1)  & ChebyKAN & FBiCGStab & 1-level & SSOR & 9-1-9 \\ \hline
Top 2 solver of $p^1$ (3d t2\_1)  & CheyKAN & FBiCGStab & 1-level & SSOR &  7-1-7 \\ \hline
Top 3 solver of $p^1$ (3d t3\_1)  & CheyKAN & FBiCGStab & 1-level & SSOR & 5-1-5 \\ \hline
%Accurate-4($\frac{1}{128}$)& ChebyKAN & SSOR & 1-level & 0 & 0.015625 \\ \hline
%Accurate-5($\frac{1}{128}$)& DeepONet & BiCGStab & 1-level & 8 & 0.5 \\  \hline
\end{tabular}
}
    \end{minipage}

    \vspace{1em} % Space between subtables

    % Second subtable
    \begin{minipage}{\linewidth}
        \centering
        \captionof{subtable}{Performance of the top-3 solvers.}
        \resizebox{\linewidth}{!}{ % Resize to fit the line width
          \begin{tabular}{c|c|c|c|c|c|c|c}
\hline \hline 
 &Error & Com. time &  \# of ite. & Conv. rate & Memory & MACs & Training time \\ 
\hline \hline
3d t1\_1& 1.56E-13 & 0.131 & 7 & 4.260 & 11886.9 & 5.830E+8 & 18157.3 \\ \hline
3d t2\_1& 1.12E-13 & 0.122 & 8 & 3.764 & 11887.3 & 5.557E+8 & 18157.3 \\ \hline
3d t3\_1&  1.74E-14 & 0.117 & 10 & 3.135 & 11898 & 5.566E+8 & 18157.3 \\ \hline
%Accurate-4($\frac{1}{128}$)& 0.0971029 & 1.24E-14 & 6 & 2.86298 & 235.11 & 178858 & 99.2164 \\ \hline
%Accurate-5($\frac{1}{128}$)& 0.0518639 & 1.32E-14 & 13 & 1.47149 & 236.639 & 178858 & 326.378 \\  \hline
\end{tabular}
        }
    \end{minipage}
    \label{top3_p1_3d_kry}
\end{table}

\begin{table}
    \centering
    \caption{Top-3 solvers by preference function $p^2$ in 3-D for Krylov-based methods.}
    
    \begin{minipage}{\linewidth}
        \centering
        \captionof{subtable}{The top-3 solvers.}
        \resizebox{\linewidth}{!}{
        \begin{tabular}{c|c|c|c|c|c}
\hline \hline 
 &Neural operator & Classical solver & Multi-grid & Relaxation & Strategies \\ 
\hline \hline
Top 1 solver of $p^1$ (3d t1\_1)  & ChebyKAN & FBiCGStab & 3-level & SSOR & 1-1-1 \\ \hline
Top 2 solver of $p^1$ (3d t2\_1) & DeepONet & FBiCGStab & 3-level & SSOR &  1-1-1 \\ \hline
Top 3 solver of $p^1$ (3d t3\_1) & JacobiKAN & FBiCGStab & 3-level & SSOR & 1-1-1 \\ \hline
%Accurate-4($\frac{1}{128}$)& ChebyKAN & SSOR & 1-level & 0 & 0.015625 \\ \hline
%Accurate-5($\frac{1}{128}$)& DeepONet & BiCGStab & 1-level & 8 & 0.5 \\  \hline
\end{tabular}
}
    \end{minipage}

    \vspace{1em} % Space between subtables

    % Second subtable
    \begin{minipage}{\linewidth}
        \centering
        \captionof{subtable}{Performance of the top-3 solvers.}
        \resizebox{\linewidth}{!}{ % Resize to fit the line width
          \begin{tabular}{c|c|c|c|c|c|c|c}
\hline \hline 
 &Error & Com. time &  \# of ite. & Conv. rate & Memory & MACs & Training time \\ 
\hline \hline
3d t1\_1& 1.79E-14 & 0.0866 & 4 & 8.047 & 13663.5 & 3.290E+8 & 18157.3\\ \hline
3d t2\_1& 1.74E-14 & 0.0906 & 4 & 8.065 & 13662 & 4.667E+9 & 8625.9 \\ \hline
3d t3\_1& 1.34E-14 & 0.0895 & 4 & 8.085 & 13664 & 3.290E+8 & 28533.4 \\ \hline
%Accurate-4($\frac{1}{128}$)& 0.0971029 & 1.24E-14 & 6 & 2.86298 & 235.11 & 178858 & 99.2164 \\ \hline
%Accurate-5($\frac{1}{128}$)& 0.0518639 & 1.32E-14 & 13 & 1.47149 & 236.639 & 178858 & 326.378 \\  \hline
\end{tabular}
        }
    \end{minipage}
    \label{top3_p2_3d_kry}
\end{table}

We implement the method introduced in~\Cref{sec-red} for re-discovering particular solvers. 
Recall that given a solver and its performance, the weights in the preference function are computed by solving the problem~\eqref{lp_discover}. 
In the following, we present the weights in the order s.t. $(\lambda_1,\lambda_2,\dots,\lambda_7)$ for ( $\lambda_1$: relative error, $\lambda_2$: computational time, $\lambda_3$: number of iterations, $\lambda_4$: convergence rate, $\lambda_5$: memory allocation, $\lambda_6$: MACs, $\lambda_7$: training time), for some particular solvers. 
The following represents the rediscovery of optimal solvers in 3-d case: 

\begin{subequations}\label{redis_kry}
\begin{equation}\label{redis_kry_a}
\left\{
\begin{aligned}
&\text{Preference function $p(\lambda_a;r) = (0, 0.012, 0, 0.007, 0.044, 0.368, 0.568)^Tr$ } \\
&\text{Optimal solver}: x^a = (\text{ChebyKAN}, \text{FCG}, \text{1-level}, \text{SSOR}, 5-1-5) \ .
\end{aligned}
\right.
\end{equation}

\begin{equation}
\left\{
\begin{aligned}
&\text{Preference function $p(\lambda_a;r) = (0, 0.406, 0, 0.281, 0.183, 0, 0.129)^Tr$ } \\
&\text{Optimal solver}: x^a = (\text{DeepONet}, \text{FGMRES}, \text{1-level}, \text{SSOR}, 5-1-5) \ .
\end{aligned}
\right.
\end{equation}

%\begin{equation}
%\left\{
%\begin{aligned}
%&\text{Preference function $p(\lambda_a;r) = (0.213, 0, 0.670, 0, 0.117, 0, 0)^Tr$ } \\
%&\text{Optimal solver}: x^a = (\text{U-Net}, \text{FBiCGStab}, \text{2-level}, \text{SSOR}, 5-1-5) \ .
%\end{aligned}
%\right.
%\end{equation}

\end{subequations}

Same as the relaxation-based methods, when the optimal solver changes, the weights in the preference function used for discovery also change significantly. Taking for instance the optimal solver $x^a$ re-discovered in~\eqref{redis_kry_a}, when the users care more about MACs and training time (with weights 0.368 and 0.568 resp. ) while assigning minimal importance to the other criteria (weights $<$ 0.05), the meta-solver ($x_1$=ChebyKAN, $x_2$=FCG, $x_3$=1-level, $x_4$=SSOR,$x_5=5-1-5$) is the optimal one. 

%% file: conclusion.tex
\section{Conclusion}\label{sec-conclu}

We propose a novel multi-objective optimization method based on Pareto optimality to model and quantify the problem of identifying an optimal meta-solver for a given task in a specific context. Building on this structure, we introduce a preference function approach for discovering the optimal solver among Pareto optimal candidates, and we employ a linear programming methodology to select a particular solver. 
We present numerical results for solving Poisson equations in 1-d/2-d (in the Appendix), and 3-d domains (in the main text), demonstrating the effectiveness and potential of our approach. We plan to extend this automation to identify optimal solvers for nonlinear systems and other applications.

%% file: appendix.tex
\appendix
\section{Proofs of Proposition}\label{appd_imp}
\subsection{Proof of~\Cref{mono_proper}}

We prove for the set of Pareto solution, then it follows straightforward for the weak Pareto sense. 

For~\eqref{mono_l}, take a $x \notin \Pp_{\tilde{L}}(X)$, then there exists a $x' \in X$ such that $\tilde{L}_k(x')\leq \tilde{L}_k(x)$ for every $k \in \{1,\dots,N,\dots, N+\tilde{N} \}$. In particular, we have $L_k(x') \leq L_k(x)$ for every $k \in \{1,\dots,N\}$.  By~\Cref{assump1}, this implies $x' \succeq x$ in the sense of minimizing $L$ over $X$, so that $x \notin \Pp_L(X)$. We conclude that $\Pp_L(X) \subseteq \Pp_{\tilde{L}}(X)$. 

For~\eqref{mono_x1}, take a $x\notin (\Pp_L(X) \cup (\tilde{X} \setminus X))$, then $x \in X$ and $x\notin \Pp_L(X)$. By definition, there exists a $x' \in X \subseteq \tilde{X}$ such that $x' \succeq x$, thus $x \notin \Pp_L(\tilde{X})$. We conclude that $\Pp_L(\tilde{X}) \subseteq \Pp_L(X) \cup (\tilde{X} \setminus X)$. For~\eqref{mono_x2}, it is sufficient to notice that $\Pp_L(\Pp_L(\tilde{X})) = \Pp_L(\tilde{X})$ and $\Pp_L(A \cup B) = \Pp_L( \Pp_L(A) \cup B)$ for any $A,B \subseteq \tilde{X}$. 
$\qqed$
\subsection{Proof of~\Cref{propo_prefer}}

Let $\bar{x}$ be a minimizer of $g \circ L$. Assume $L(\bar{x}) \notin \F_L(X)$. Then there exists a $x'$ dominates $\bar{x}$, that is, 
\begin{equation}\forall \ i, \ L_i(x') \leq L_i(\bar{x}), \text{ and } \exists \ j \text{ such that } L_j(x') < L_j(\bar{x}) \ . \end{equation}
Since $g$ is increasing, we have $g \circ L(x') < g \circ L(\bar{x})$, so contradiction. 
$\qqed$

\section{Further numerical results}\label{app_numerics}

\subsection{Relaxation-based methods for solving 1-d, 2-d Poisson equations}

\begin{table}
    \centering
    \tiny
    \caption{ Composition of Pareto optimal solvers by counting the number of elements in each dimension, with the highest number highlighted, for relaxation-based methods.}
    \begin{minipage}{\linewidth}
        \centering
        \captionof{subtable}{Different neural operators}
        \resizebox{\linewidth}{!}{ 
        \begin{tabular}{c| c c c c c c c} 
\hline
Neural Op & DeepONet & U-Net & FNO & Transformer & KAN & JacobiKAN& ChebyKAN  \\
\hline
\# in Pareto opt in 1D  & 7 & 5 & 11 & 11 & 13 & \textbf{17} & 12 \\
\hline
\# in Pareto opt in 2D  & 16 & \textbf{18} & 4 & 5 & 15 & 17 & 0 \\
%\hline
%\# in Pareto opt in 3D  & 20 & 18 & 7 & 5 & 21 & \textbf{27} & 21 \\
\hline
%\# in Pareto opt., $\frac{1}{128}$ & 28 & 20 & 7 & 8 & 14 & 36
\end{tabular}
}
    \end{minipage}

   % \vspace{1em} % Space between subtables

    % Second subtable
    \begin{minipage}{\linewidth}
        \centering
        \captionof{subtable}{Different classical iterative solvers.}
           \begin{tabular}{c| c c c c    } 
\hline
Classical solvers & Jacobi & Gauss-Siedel & SOR & SSOR  \\
\hline
\# in Pareto opt in 1D & 2 & 19 & 22 & \textbf{33}  \\
\hline
\# in Pareto opt in 2D & 2 & 12 & 23 & \textbf{38}  \\
%\hline
%\# in Pareto opt in 3D & 8 & 14 & 37 & \textbf{60}  \\
\hline
%\# in Pareto opt., $\frac{1}{128}$ & 1 & 8 & 26 & 27 & 35 & 16 & 0 
\end{tabular}
        
    \end{minipage}

\begin{minipage}{\linewidth}
        \centering
        \captionof{subtable}{Different number of multi-grids.}
           \begin{tabular}{c| c c c  } 
\hline
Multi-grid method& NO & 2-level & 3-level \\
\hline
\# in Pareto opt in 1D & 1 & \textbf{56} & 19 \\
\hline
\# in Pareto opt in 2D & 13 & \textbf{43} & 19 \\
%\hline
%\# in Pareto opt in 3D & 15 & \textbf{63} & 41 \\
\hline
%\# in Pareto opt., $\frac{1}{128}$ & 62 & 42 & 9 \\
\end{tabular}
        
    \end{minipage}

\begin{minipage}{\linewidth}
        \centering
        \captionof{subtable}{Different proportion of neural operators.}
    \begin{tabular}{c| c c c c c c c } 
\hline
Ratio of N.O. & $\frac{1}{2}$ & $\frac{1}{4}$ & $\frac{1}{8}$ & $\frac{1}{16}$ & $\frac{1}{32}$ & $\frac{1}{64}$ & $\frac{1}{128}$  \\ 
\hline
\# in Pareto opt in 1D & 0 & 0 & 5 & 13 & {13} & {21} & \textbf{24} \\ 
\hline
\# in Pareto opt in 2D & 0 & 0 & 3 & 9 & 14 & 24 & \textbf{25} \\ 
%\hline
%\# in Pareto opt in 3D & 0 & 2 & 12 & 18 & 27 & 28 & \textbf{32} \\ 
\hline
%\# in Pareto opt., $\frac{1}{128}$ & 6 & 6 & 10 & 22 & 23 & 23 & 23 \\ 
\end{tabular}

       \end{minipage}
    \label{relax_com_app}
\end{table}

\begin{figure}
\centering
\subfigure[Computational time -- Relative error -- \# of iterations]{
\includegraphics[width=0.3\textwidth]{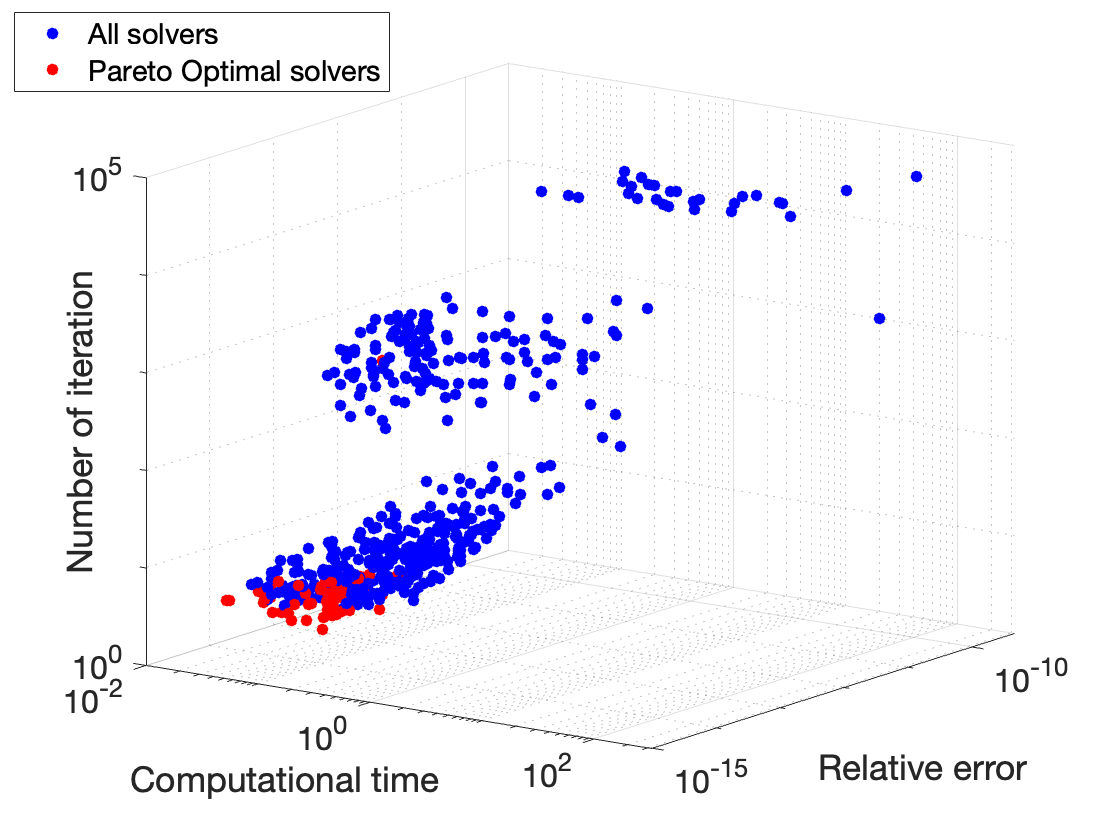}}
\subfigure[Relative error -- MACs -- Memory allocation]{
\includegraphics[width=0.3\textwidth]{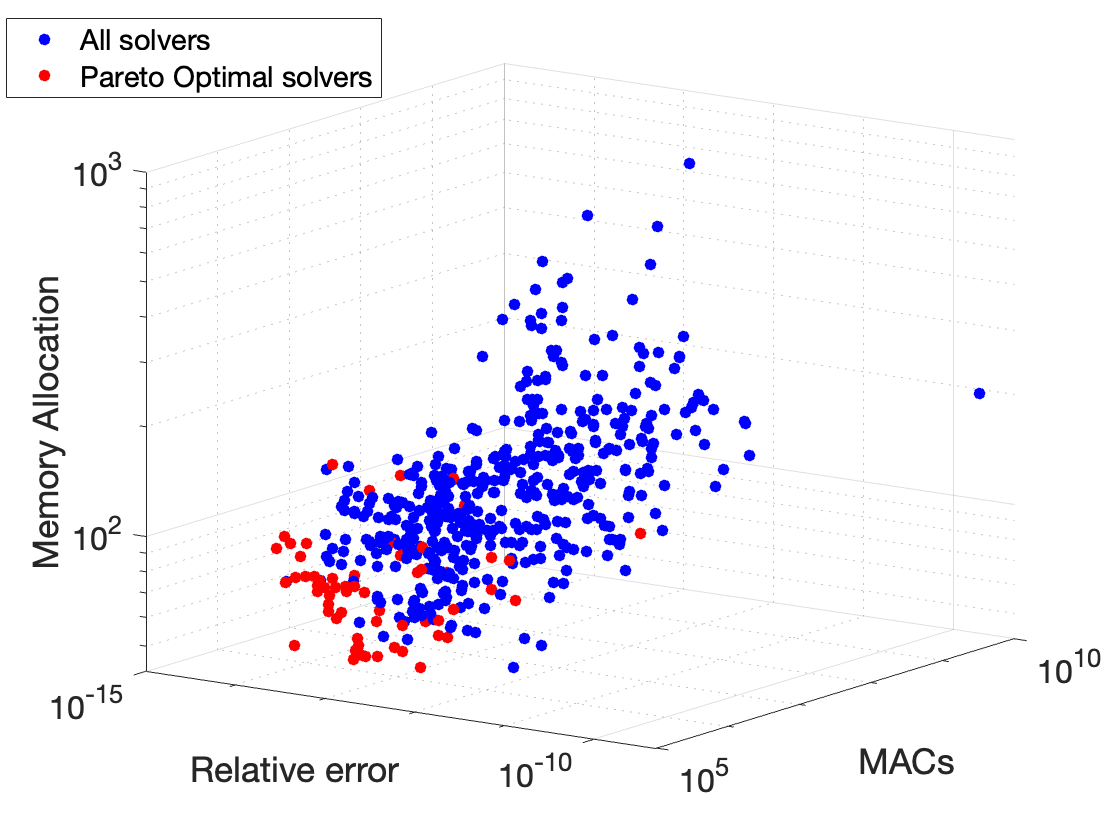}}
\subfigure[Convergence rate -- MACs -- \# of iterations]{
\includegraphics[width=0.3\textwidth]{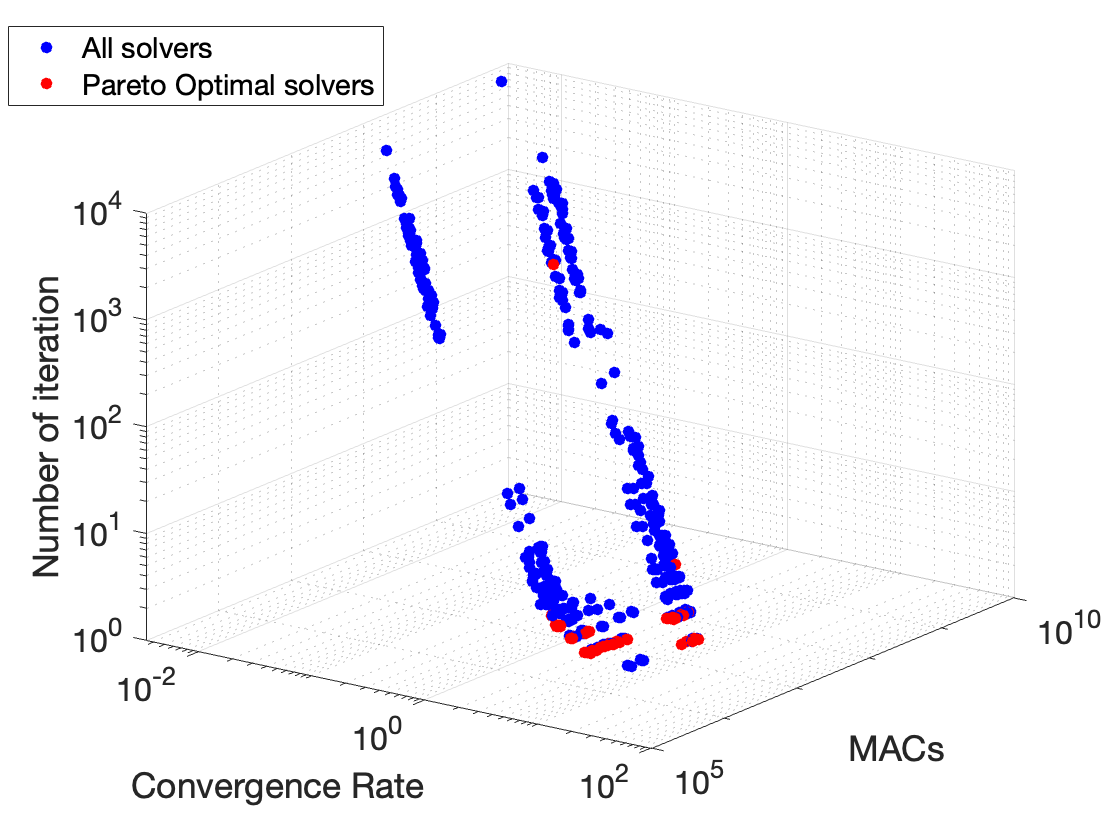}}
\caption{Projection of Pareto front onto three dimensional criteria, for relaxation-based methods solving 1-d Poisson equation.}
\label{front_tei_12}
\end{figure}

\begin{figure}
\centering
\subfigure[Computational time -- Relative error -- \# of iterations]{
\includegraphics[width=0.3\textwidth]{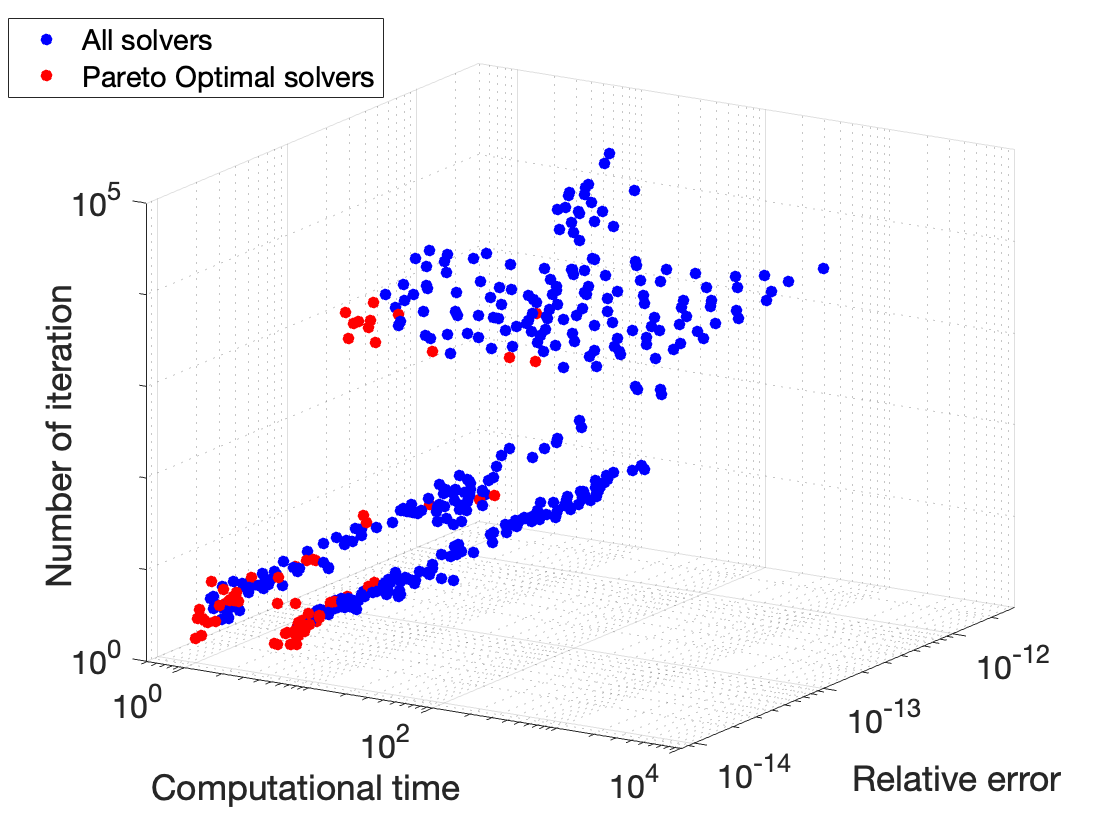}}
\subfigure[Relative error -- MACs -- Memory allocation]{
\includegraphics[width=0.3\textwidth]{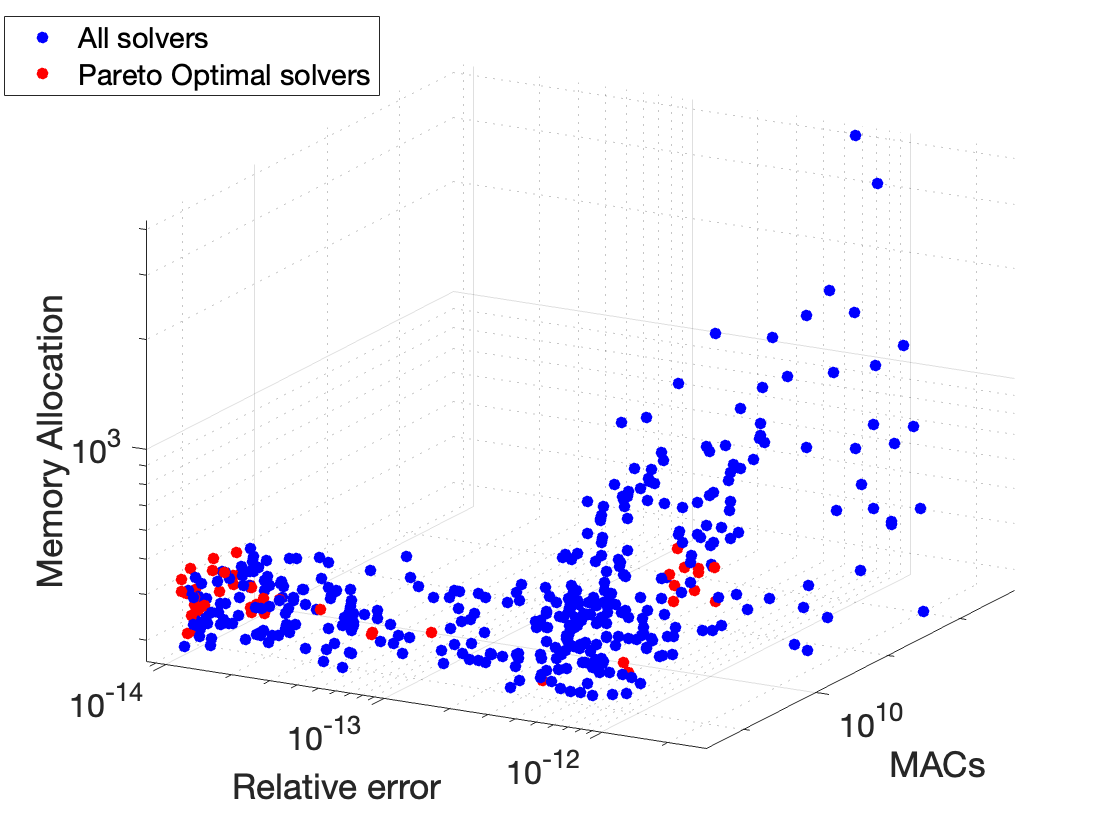}}
\subfigure[Convergence rate -- MACs -- \# of iterations]{
\includegraphics[width=0.3\textwidth]{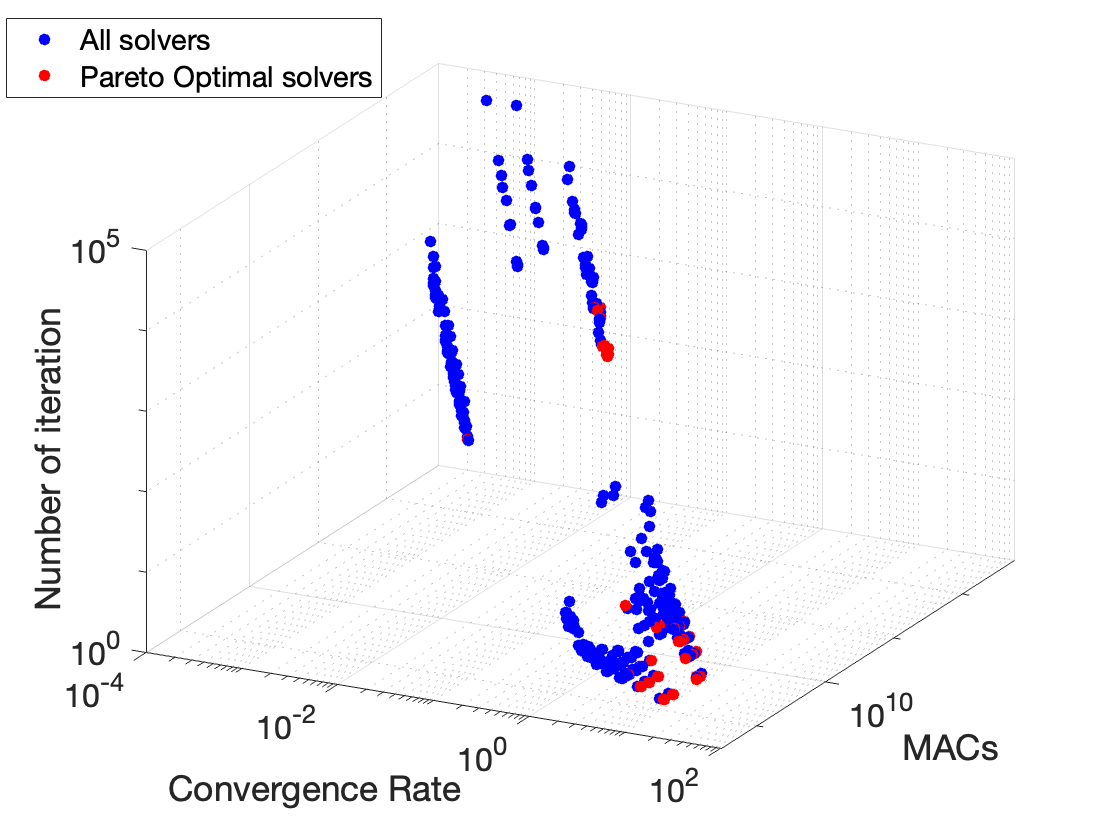}}
\caption{Projection of Pareto front onto three dimensional criteria, for relaxation-based methods solving 2-d Poisson equation.}
\label{front_emm_12}
\end{figure}

\begin{table}
\tiny
    \centering
    \caption{Top-3 solvers by preference function $p^1$ in 1-D for relaxation-based methods.}
    \begin{minipage}{\linewidth}
        \centering
        \captionof{subtable}{The top-3 solvers.}
       
        \begin{tabular}{c|c|c|c|c}
\hline \hline 
 &Neural operator & Classical solver & Multi-grid &  Ratio \\ 
\hline \hline
Top 1 solver of $p^1$ (1d t1\_1) in 1-d &Transformer &SSOR & 2-level & $\frac{1}{32}$\\ \hline
Top 2 solver of $p^1$ (1d t2\_1) in 1-d &Transformer & SOR & 2-level & $\frac{1}{128}$ \\ \hline
Top 3 solver of $p^1$ (1d t3\_1) in 1-d & Transformer & Gauss-seidel & 2-level &  $\frac{1}{128}$ \\ \hline
%Accurate-4($\frac{1}{128}$)& ChebyKAN & SSOR & 1-level & 0 & 0.015625 \\ \hline
%Accurate-5($\frac{1}{128}$)& DeepONet & BiCGStab & 1-level & 8 & 0.5 \\  \hline
\end{tabular}

    \end{minipage}

  %  \vspace{1em} % Space between subtables

    % Second subtable
    \begin{minipage}{\linewidth}
        \centering
        \captionof{subtable}{Performance of the top-3 solvers.}
        % Resize to fit the line width
           \begin{tabular}{c|c|c|c|c|c|c|c}
\hline \hline 
 &Error & Com. time &  \# of ite. & Conv. rate & Memory & MACs & Training time \\ 
\hline \hline
1d t1\_1& 2.45E-14 & 0.0473 & 3 & 10.763 & 49.762 & 6.85E+6 & 23.669\\ \hline
1d t2\_1& 2.46E-14 & 0.0484 & 3 & 10.889 & 50.856 & 7.35E+6 & 23.669 \\ \hline
1d t3\_1& 2.16E-14 & 0.0484 & 3 & 10.853 & 51.802 & 7.35E+6 & 23.669 \\ \hline
%Accurate-4($\frac{1}{128}$)& 0.0971029 & 1.24E-14 & 6 & 2.86298 & 235.11 & 178858 & 99.2164 \\ \hline
%Accurate-5($\frac{1}{128}$)& 0.0518639 & 1.32E-14 & 13 & 1.47149 & 236.639 & 178858 & 326.378 \\  \hline
\end{tabular}
        
    \end{minipage}
    \label{top3_p1_1d}
\end{table}

\begin{table}
\tiny
    \centering
    \caption{Top-3 solvers by preference function $p^1$ in 2-D for relaxation-based methods.}
    
    \begin{minipage}{\linewidth}
        \centering
        \captionof{subtable}{The top-3 solvers.}
        \resizebox{\linewidth}{!}{
        \begin{tabular}{c|c|c|c|c}
\hline \hline 
 &Neural operator & Classical solver & Multi-grid &  Ratio \\ 
\hline \hline
Top 1 solver of $p^1$ (2d t1\_1) in 2-d & DeepONet & SSOR & 2-level & $\frac{1}{128}$ \\ \hline
Top 2 solver of $p^1$ (2d t2\_1) in 2-d & U-Net & SSOR & 3-level &  $\frac{1}{128}$ \\ \hline
Top 3 solver of $p^1$ (2d t3\_1) in 2-d & U-Net & SSOR & 2-level &  $\frac{1}{128}$ \\ \hline
%Accurate-4($\frac{1}{128}$)& ChebyKAN & SSOR & 1-level & 0 & 0.015625 \\ \hline
%Accurate-5($\frac{1}{128}$)& DeepONet & BiCGStab & 1-level & 8 & 0.5 \\  \hline
\end{tabular}
}
    \end{minipage}

   % \vspace{1em} % Space between subtables

    % Second subtable
    \begin{minipage}{\linewidth}
        \centering
        \captionof{subtable}{Performance of the top-3 solvers.}
        % Resize to fit the line width
          \begin{tabular}{c|c|c|c|c|c|c|c}
\hline \hline 
 &Error & Com. time &  \# of ite. & Conv. rate & Memory & MACs & Training time \\ 
\hline \hline
2d t1\_1& 1.01E-14 & 1.071 & 2 & 16.878 & 351.142 & 1.650E+9 & 6291.1 \\ \hline
2d t2\_1& 4.28E-14 & 1.373 & 2 & 15.382 & 337.712 & 2.168E+9 & 2207.45 \\ \hline
2d t3\_1& 3.87E-14 & 1.081 & 2 & 15.416 & 346.537 & 1.827E+9 & 2207.45 \\ \hline
%Accurate-4($\frac{1}{128}$)& 0.0971029 & 1.24E-14 & 6 & 2.86298 & 235.11 & 178858 & 99.2164 \\ \hline
%Accurate-5($\frac{1}{128}$)& 0.0518639 & 1.32E-14 & 13 & 1.47149 & 236.639 & 178858 & 326.378 \\  \hline
\end{tabular}
        
    \end{minipage}
    \label{top3_p1_2d}
\end{table}

\begin{figure}
\centering
\subfigure[1-d case by $p^1$]{
\includegraphics[width=0.23\textwidth]{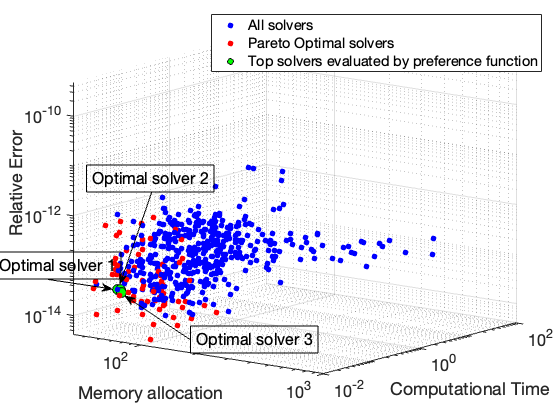}}
\subfigure[2-d case by $p^1$]{
\includegraphics[width=0.23\textwidth]{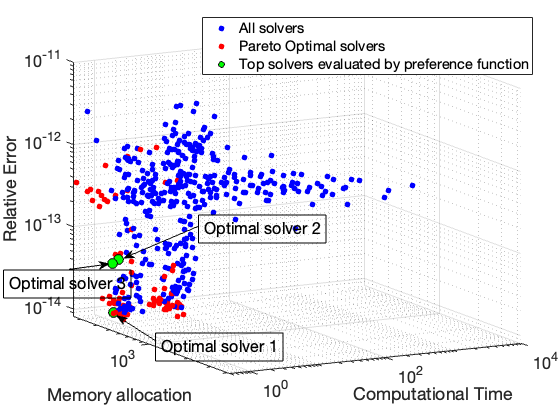}}
\subfigure[1-d case by $p^2$]{
\includegraphics[width=0.23\textwidth]{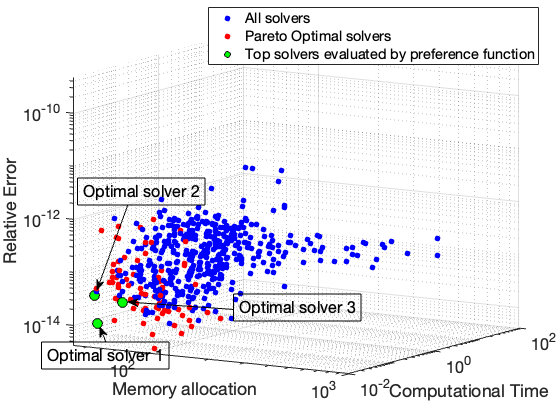}}
\subfigure[2-d case by $p^2$]{
\includegraphics[width=0.23\textwidth]{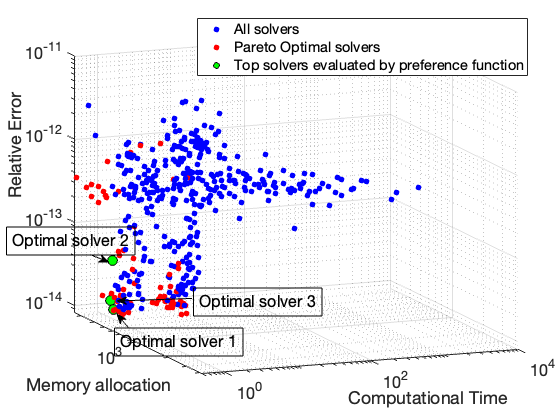}}
%\subfigure[Solving 3d Poisson Equation]{
%\includegraphics[width=0.3\textwidth]{figs/p1_relax_3d.png}}
\caption{The top-3 optimal solvers discovered by preference function $p^1$ (left 2), and $p^2$ (right 2) in the three-dimensional projection: Computational Time -- Relative Error -- Memory allocation of Pareto front, for relaxation-based methods.}
\label{front_p1_relax}
\end{figure}

\begin{table}
\tiny
    \centering
    \caption{Top-3 solvers by preference function $p^2$ in 1-D for relaxation-based method.}
    \begin{minipage}{\linewidth}
        \centering
        \captionof{subtable}{The top-3 solvers.}
        
        \begin{tabular}{c|c|c|c|c}
\hline \hline 
 &Neural operator & Classical solver & Multi-grid &  Ratio \\ 
\hline \hline
Top 1 solver of $p^2$ (t1\_2 1d) in 1-d & DeepONet & SOR & 2-level & $\frac{1}{128}$ \\ \hline
Top 2 solver of $p^2$ (t2\_2 1d) in 1-d& DeepONet & SOR & 2-level &  $\frac{1}{64}$ \\ \hline
Top 3 solver of $p^2$ (t3\_2 1d) in 1-d& DeepONet & Gauss-Seidel & 2-level &  $\frac{1}{128}$ \\ \hline
%Accurate-4($\frac{1}{128}$)& ChebyKAN & SSOR & 1-level & 0 & 0.015625 \\ \hline
%Accurate-5($\frac{1}{128}$)& DeepONet & BiCGStab & 1-level & 8 & 0.5 \\  \hline
\end{tabular}

    \end{minipage}

  %  \vspace{1em} % Space between subtables

    % Second subtable
    \begin{minipage}{\linewidth}
        \centering
        \captionof{subtable}{Performance and rank of performance of the top-3 solvers.}
        \resizebox{\linewidth}{!}{ % Resize to fit the line width
           \begin{tabular}{c|c|c|c|c|c|c|c}
\hline \hline 
 &Error & Com. time &  \# of ite. & Conv. rate & Memory & MACs & Training time \\ 
\hline \hline
t1\_2 1d & 1.07E-14 & 0.0143 & 3 & 10.679 & 52.063 & 9.41E+6 & 174.409 \\ \hline
t2\_2 1d & 3.59E-14 & 0.0137 & 3 & 10.771 & 51.067 & 9.02E+6 & 174.409 \\ \hline
t3\_2 1d &  3.10E-14 & 0.0143 & 3 & 10.363 & 70.133 & 9.41E+6 & 174.409 \\ \hline
%Accurate-4($\frac{1}{128}$)& 0.0971029 & 1.24E-14 & 6 & 2.86298 & 235.11 & 178858 & 99.2164 \\ \hline
%Accurate-5($\frac{1}{128}$)& 0.0518639 & 1.32E-14 & 13 & 1.47149 & 236.639 & 178858 & 326.378 \\  \hline
\end{tabular}
        }
    \end{minipage}
    \label{top3_p2_1d}
\end{table}

\begin{table}
\tiny
    \centering
    \caption{Top-3 solvers by preference function $p^2$ in 2-D for relaxation-based method.}
    \begin{minipage}{\linewidth}
        \centering
        \captionof{subtable}{The top-3 solvers.}
        
        \begin{tabular}{c|c|c|c|c}
\hline \hline 
 &Neural operator & Classical solver & Multi-grid &  Ratio \\ 
\hline \hline
Top 1 solver of $p^2$ (t1\_2 2d) in 2-d & DeepONet & SSOR & 2-level & $\frac{1}{128}$ \\ \hline
Top 2 solver of $p^2$ (t2\_2 2d) in 2-d& U-Net & SSOR & 2-level &  $\frac{1}{128}$ \\ \hline
Top 3 solver of $p^2$ (t3\_2 2d) in 2-d& DeepONet & Gauss-Seidel & 2-level &  $\frac{1}{128}$ \\ \hline
%Accurate-4($\frac{1}{128}$)& ChebyKAN & SSOR & 1-level & 0 & 0.015625 \\ \hline
%Accurate-5($\frac{1}{128}$)& DeepONet & BiCGStab & 1-level & 8 & 0.5 \\  \hline
\end{tabular}

    \end{minipage}

   % \vspace{1em} % Space between subtables

    % Second subtable
    \begin{minipage}{\linewidth}
        \centering
        \captionof{subtable}{Performance and rank of performance of the top-3 solvers.}
        % Resize to fit the line width
           \begin{tabular}{c|c|c|c|c|c|c|c}
\hline \hline 
 &Error & Com. time &  \# of ite. & Conv. rate & Memory & MACs & Training time \\ 
\hline \hline
t1\_2 2d & 1.01E-14 & 1.071 & 2 & 16.878 & 351.142 & 1.650E+9 & 6291.1 \\ \hline
t2\_2 2d & 3.87E-14 & 1.081 & 2 & 15.416 & 346.537 & 1.827E+9 & 2207.5 \\ \hline
t3\_2 2d & 1.17E-14 & 1.233 & 4 & 8.223 & 304.189 & 2.873E+9 & 6291.1 \\ \hline
%Accurate-4($\frac{1}{128}$)& 0.0971029 & 1.24E-14 & 6 & 2.86298 & 235.11 & 178858 & 99.2164 \\ \hline
%Accurate-5($\frac{1}{128}$)& 0.0518639 & 1.32E-14 & 13 & 1.47149 & 236.639 & 178858 & 326.378 \\  \hline
\end{tabular}
        
    \end{minipage}
    \label{top3_p2_2d}
\end{table}

\eqref{redis_relax_1d} presents the rediscovery of optimal solvers in 1-d for relaxation-based methods, and \eqref{redis_relax_2d} presents the rediscovery of optimal solvers in 2-d.  
\begin{subequations}\label{redis_relax_1d}
\begin{equation}
\left\{
\begin{aligned}
&\text{Preference function $p(\lambda_a;r) = (0.95, 0.01, 0, 0.04, 0, 0,  0)^Tr$ } \\
&\text{Optimal solver}: x^a = (\text{U-Net}, \text{Jaccobi}, \text{2-level}, \frac{1}{128}) \ .
\end{aligned}
\right.
\end{equation}

\begin{equation}
\left\{
\begin{aligned}
&\text{Preference function $p(\lambda_b;r) = (0.01, 0.02, 0, 0.36, 0.59, 0.02,  0.24)^Tr$ } \\
&\text{Optimal solver}: x^b = (\text{KAN}, \text{SOR}, \text{2-level}, \frac{1}{128}) \ .
\end{aligned}
\right.
\end{equation}

\end{subequations}

\begin{subequations}\label{redis_relax_2d}
\begin{equation}
\left\{
\begin{aligned}
&\text{Preference function $p(\lambda_a;r) = (0.013, 0, 0, 0, 0, 0.553,  0.435)^Tr$ } \\
&\text{Optimal solver}: x^a = (\text{JacobiKAN}, \text{SOR}, \text{2-level}, \frac{1}{16}) \ .
\end{aligned}
\right.
\end{equation}

\begin{equation}
\left\{
\begin{aligned}
&\text{Preference function $p(\lambda_b;r) = (0, 0.65, 0, 0.05, 0.22, 0.08,  0)^Tr$ } \\
&\text{Optimal solver}: x^b = (\text{U-Net}, \text{Gauss-Seidel}, \text{2-level}, \frac{1}{32}) \ .
\end{aligned}
\right.
\end{equation}

\end{subequations}

\subsection{Krylov-based methods for solving 1-d, 2-d Poisson equations}

\begin{table}
    \centering
    \tiny
    \caption{Composition of Pareto optimal solvers by counting the number of elements in each dimension, with the highest number highlighted, for Krylov-based methods.}
    \begin{minipage}{\linewidth}
        \centering
        \captionof{subtable}{Different neural operators.}
        
        \begin{tabular}{c| c c c  c c } 
\hline
Neural Op & DeepONet & U-Net  & KAN & JacobiKAN& ChebyKAN \\
\hline
\# in Pareto opt in 1-d & 27 & 5 & 11 & \textbf{29} & 8 \\
\hline
\# in Pareto opt in 2-d & 25 & 42 & 56 & \textbf{62} & 36 \\
%\hline
%\# in Pareto opt in 3-d & \textbf{61} & 31 & 39 & 28 & 51 \\
\hline
%\# in Pareto opt., $\frac{1}{128}$ & 28 & 20 & 7 & 8 & 14 & 36
\end{tabular}

    \end{minipage}

  %  \vspace{1em} % Space between subtables

    % Second subtable
    \begin{minipage}{\linewidth}
        \centering
        \captionof{subtable}{Different Krylov solvers.}
           \begin{tabular}{c| c c c c    } 
\hline
Classical solvers & FGMRES & FCG & FBiCGStab \\
\hline
\# in Pareto opt in 1-d & 22 & 13 & \textbf{45}  \\
\hline
\# in Pareto opt in 2-d & 53 & 50 & \textbf{118}  \\
%\hline
%\# in Pareto opt in 3-d & 61 & 39 & \textbf{110}  \\
\hline
%\# in Pareto opt., $\frac{1}{128}$ & 1 & 8 & 26 & 27 & 35 & 16 & 0 
\end{tabular}
        
    \end{minipage}

    \begin{minipage}{\linewidth}
        \centering
        \captionof{subtable}{Different smoothers.}
           \begin{tabular}{c| c c c c    } 
\hline
Smoother & GS & Jacobi & SOR & SSOR  \\
\hline
\# in Pareto opt in 1-d & 15 & 15 & 15 & \textbf{35}  \\
\hline
\# in Pareto opt in 2-d & 43 & 39 & 56 & \textbf{83}  \\
%\hline
%\# in Pareto opt in 3-d & 16 & 7 & 49 & \textbf{138}  \\
\hline
%\# in Pareto opt., $\frac{1}{128}$ & 1 & 8 & 26 & 27 & 35 & 16 & 0 
\end{tabular}
        
    \end{minipage}

\begin{minipage}{\linewidth}
        \centering
        \captionof{subtable}{Different strategies of applying smoothers.}
          \begin{tabular}{c| c c c c c    } 
\hline
Classical solvers & 1-1-1 & 3-1-3 & 5-1-5 & 7-1-7 & 9-1-9  \\
\hline
\# in Pareto opt in 1-d & 20 & 6 & 10 & 19 & \textbf{25} \\
\hline
\# in Pareto opt in 2-d & 16 & \textbf{57} & 41 & 50 & \textbf{57} \\
%\hline
%\# in Pareto opt in 3-d & 34 & 43 & 44 & \textbf{54} & 35 \\
\hline
%\# in Pareto opt., $\frac{1}{128}$ & 1 & 8 & 26 & 27 & 35 & 16 & 0 
\end{tabular}
        
    \end{minipage}

\begin{minipage}{\linewidth}
        \centering
        \captionof{subtable}{Different levels of multigrid.}
    \begin{tabular}{c| c c c     } 
\hline
Levels in multi-grid & 1-level & 2-level &  3-level   \\
\hline
\# in Pareto opt in 1-d & \textbf{46} & 28 & 6  \\
\hline
\# in Pareto opt in 2-d & 68 & \textbf{88} & 65  \\
%\hline
%\# in Pareto opt in 3-d & 57 & \textbf{83} & 70  \\
\hline
%\# in Pareto opt., $\frac{1}{128}$ & 1 & 8 & 26 & 27 & 35 & 16 & 0 
\end{tabular}

       \end{minipage}
    \label{kry_com_app}
\end{table}

\begin{figure}
\centering
\subfigure[Computational time -- Relative error -- \# of iterations]{
\includegraphics[width=0.3\textwidth]{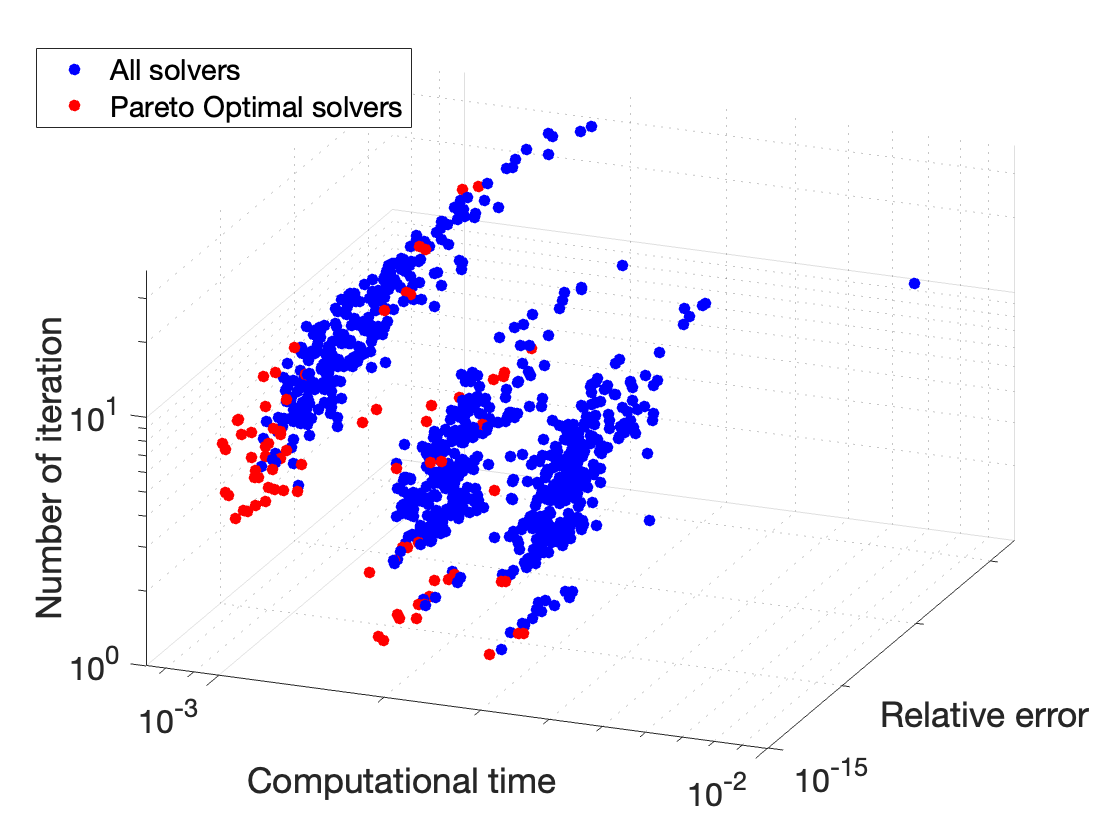}}
\subfigure[Relative error -- MACs -- Memory allocation]{
\includegraphics[width=0.3\textwidth]{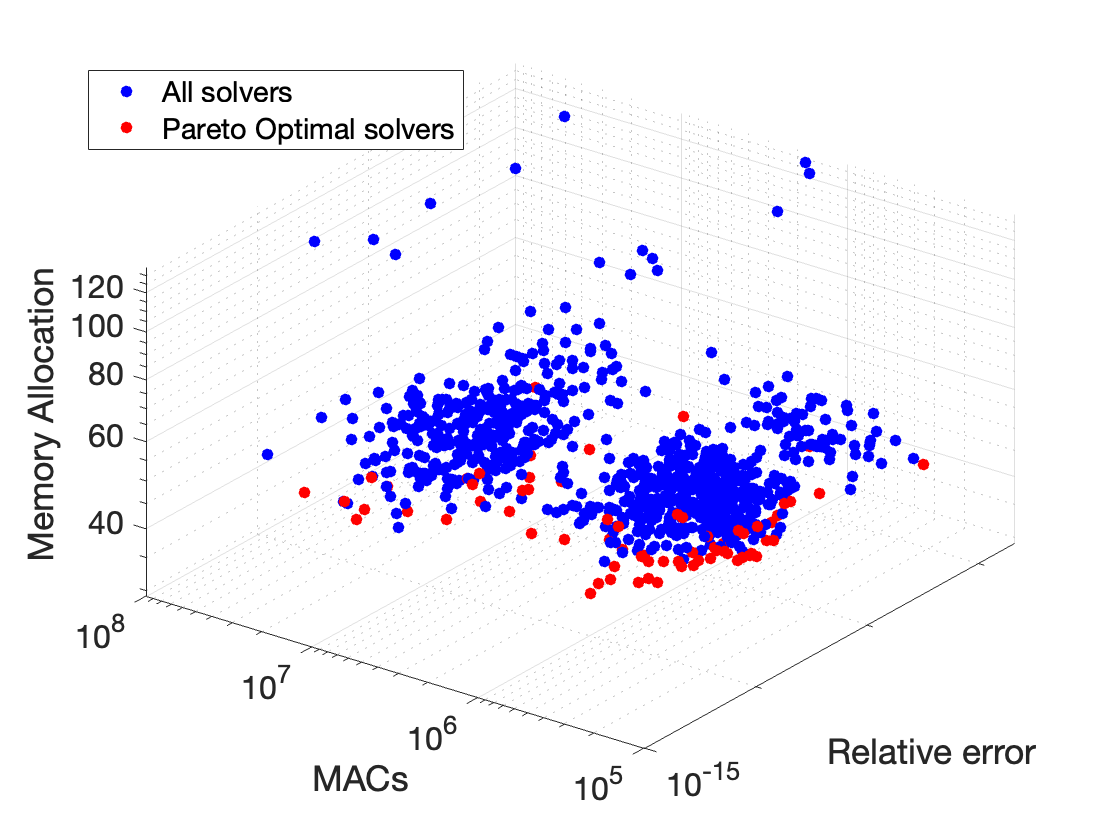}}
\subfigure[Convergence rate -- MACs -- \# of iterations]{
\includegraphics[width=0.3\textwidth]{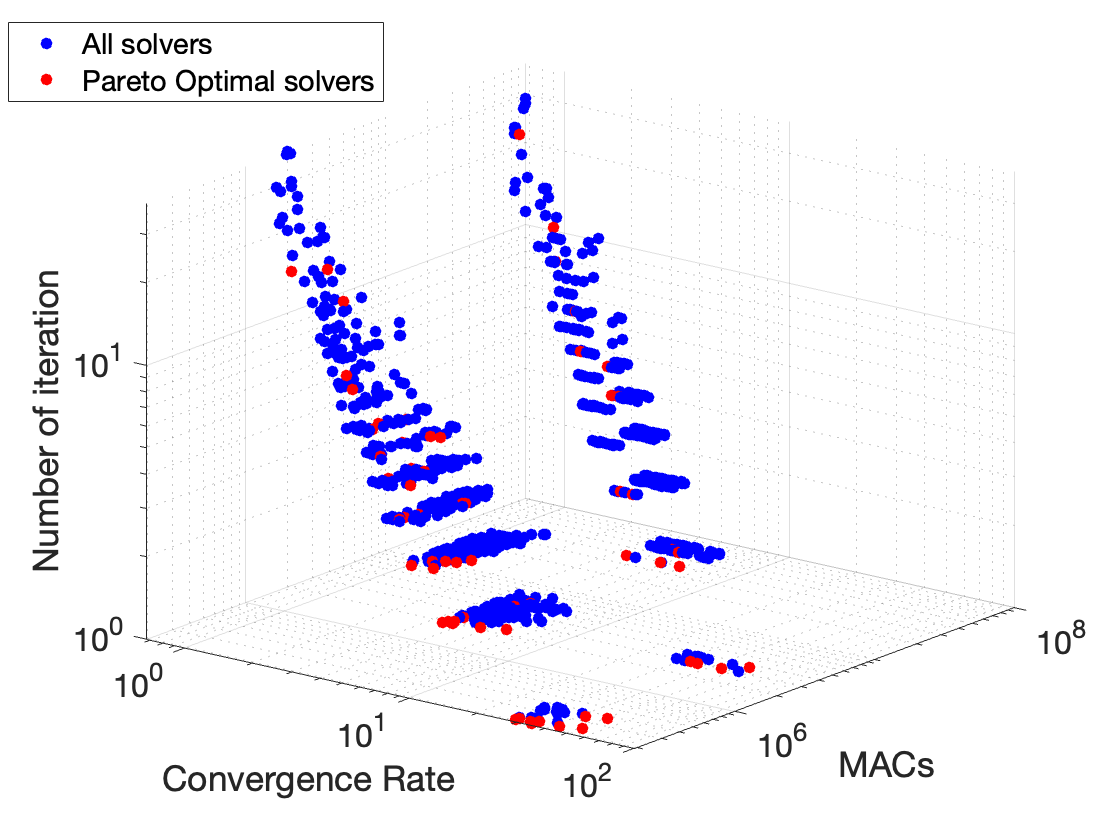}}
\caption{Projection of Pareto front onto three dimensional criteria, for Krylov-based methods, solving 1-d Poisson equation.}
\label{front_tei_12_kry}
\end{figure}

\begin{figure}
\centering
\subfigure[Computational time -- Relative error -- \# of iterations]{
\includegraphics[width=0.3\textwidth]{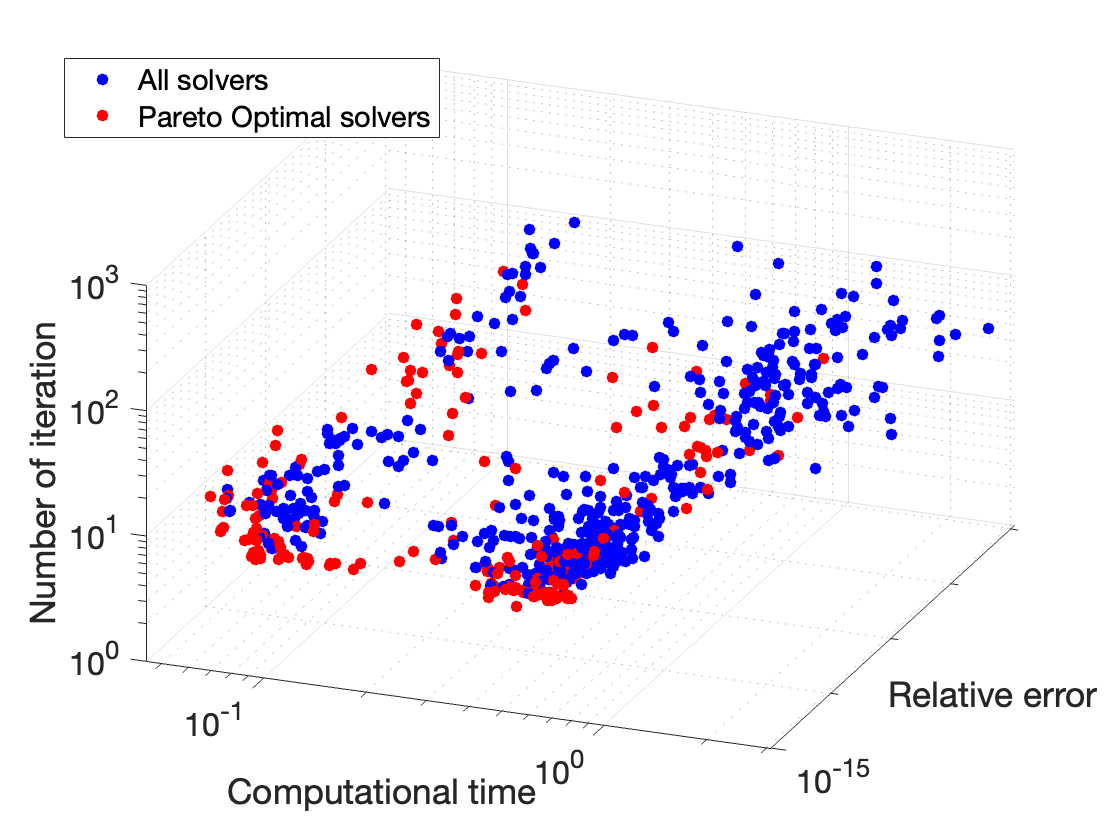}}
\subfigure[Relative error -- MACs -- Memory allocation]{
\includegraphics[width=0.3\textwidth]{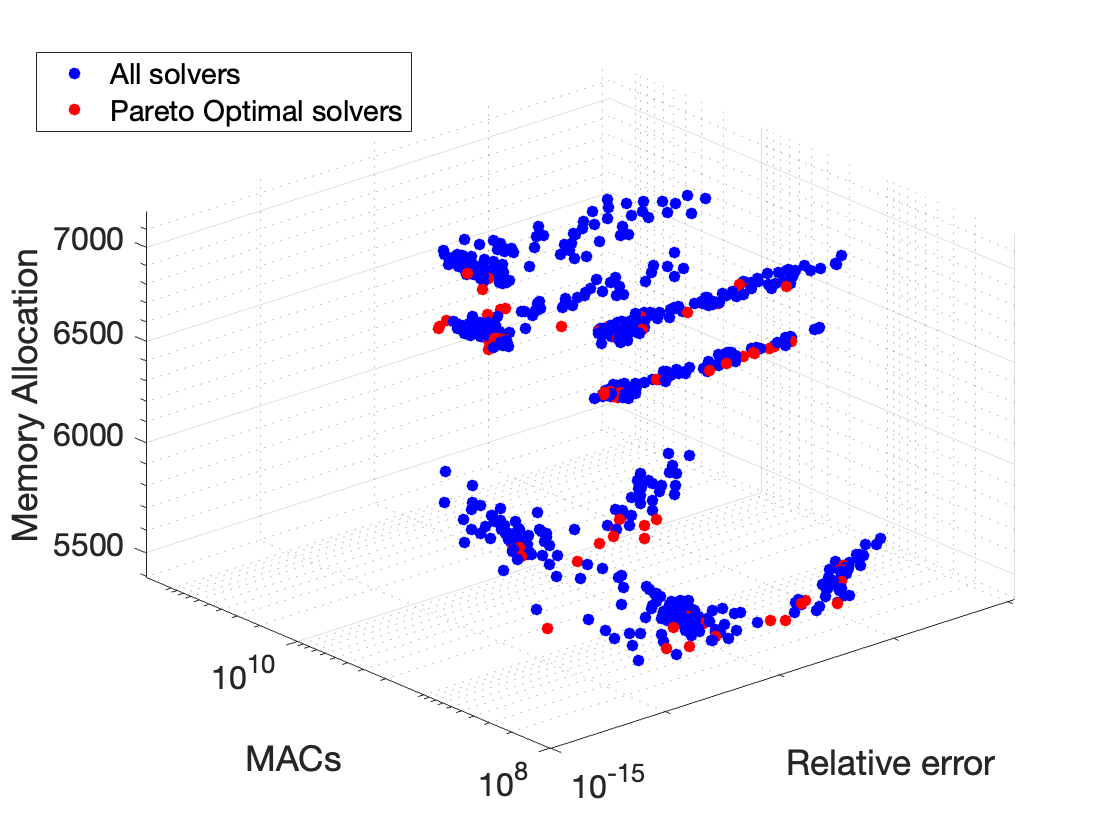}}
\subfigure[Convergence rate -- MACs -- \# of iterations]{
\includegraphics[width=0.3\textwidth]{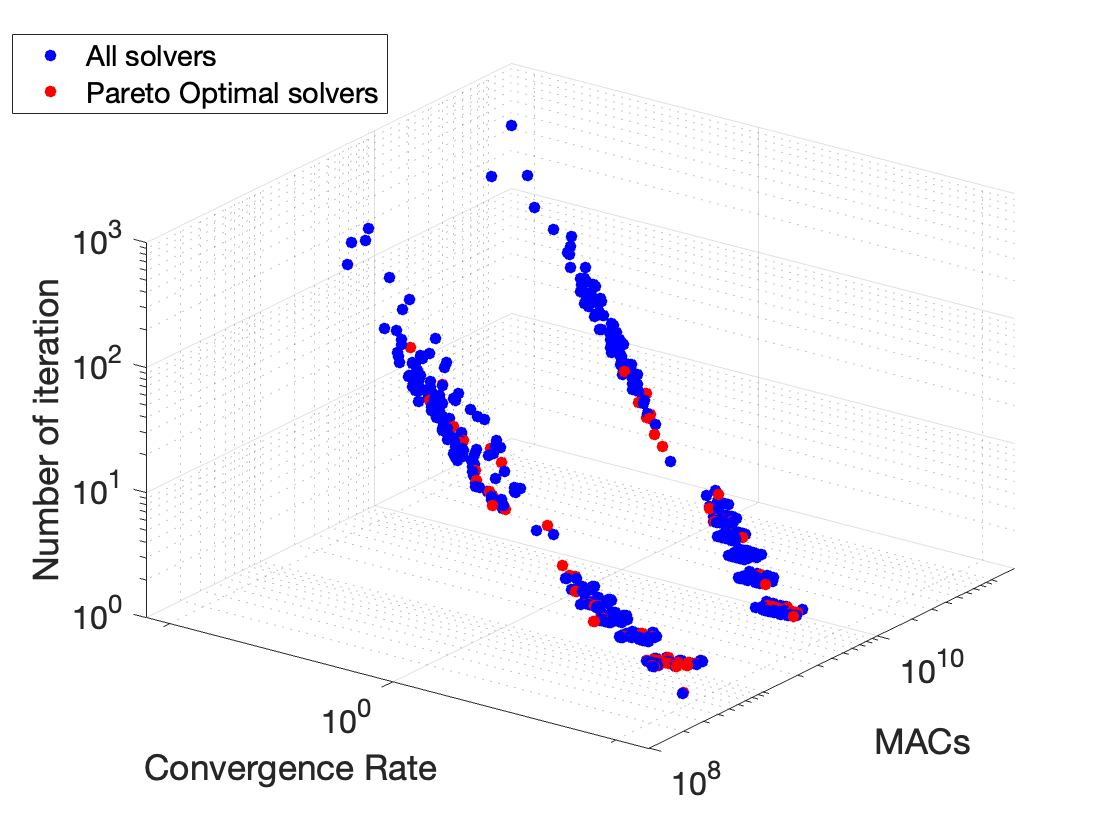}}
\caption{Projection of Pareto front onto three dimensional criteria, for Krylov-based methods, solving 2-d Poisson equation.}
\label{front_emm_12_kry}
\end{figure}

\begin{table}
\tiny
    \centering
    \caption{Top-3 solvers by preference function $p^1$ in 1-D for Krylov-based methods.}
    \begin{minipage}{\linewidth}
        \centering
        \captionof{subtable}{The top-3 solvers.}
       
        \begin{tabular}{c|c|c|c|c|c}
\hline \hline 
 &Neural operator & Classical solver & Multi-grid & Relaxation & Strategies \\ 
\hline \hline
Top 1 solver of $p^1$ (1d t1\_1) in 1-d & JacobiKAN & FBiCGStab & 1-level & SSOR & 9-1-9 \\ \hline
Top 2 solver of $p^1$ (1d t2\_1) in 1-d & JacobiKAN & FBiCGStab & 1-level & SSOR &  7-1-7 \\ \hline
Top 3 solver of $p^1$ (1d t3\_1) in 1-d & JacobiKAN & FBiCGStab & 2-level & SSOR & 9-1-9 \\ \hline
%Accurate-4($\frac{1}{128}$)& ChebyKAN & SSOR & 1-level & 0 & 0.015625 \\ \hline
%Accurate-5($\frac{1}{128}$)& DeepONet & BiCGStab & 1-level & 8 & 0.5 \\  \hline
\end{tabular}

    \end{minipage}

 %   \vspace{1em} % Space between subtables

    % Second subtable
    \begin{minipage}{\linewidth}
        \centering
        \captionof{subtable}{Performance of the top-3 solvers.}
        \resizebox{\linewidth}{!}{ % Resize to fit the line width
           \begin{tabular}{c|c|c|c|c|c|c|c}
\hline \hline 
 &Error & Com. time &  \# of ite. & Conv. rate & Memory & MACs & Training time \\ 
\hline \hline
1d t1\_1& 4.56E-14 & 8.308E-4 & 2 & 23.650 & 43.6436 & 237228 & 56.6184 \\ \hline
1d t2\_1& 2.84E-14 & 7.749E-4 & 2 & 19.189 & 43.7031 & 218748 & 56.6184 \\ \hline
1d t3\_1&  1.40E-14 & 1.605E-3 & 1 & 50.832 & 51.4268 & 242013 & 56.6184 \\ \hline
%Accurate-4($\frac{1}{128}$)& 0.0971029 & 1.24E-14 & 6 & 2.86298 & 235.11 & 178858 & 99.2164 \\ \hline
%Accurate-5($\frac{1}{128}$)& 0.0518639 & 1.32E-14 & 13 & 1.47149 & 236.639 & 178858 & 326.378 \\  \hline
\end{tabular}
        }
    \end{minipage}
    \label{top3_p1_1d_kry}
\end{table}

\begin{table}
\tiny
    \centering
    \caption{Top-3 solvers by preference function $p^1$ in 2-D for Krylov-based methods.}
    
    \begin{minipage}{\linewidth}
        \centering
        \captionof{subtable}{The top-3 solvers.}
        
       \begin{tabular}{c|c|c|c|c|c}
\hline \hline 
 &Neural operator & Classical solver & Multi-grid & Relaxation & Strategies \\ 
\hline \hline
Top 1 solver of $p^1$ (2d t1\_1) in 2-d & U-Net & FBiCGStab & 2-level & SSOR & 7-1-7 \\ \hline
Top 2 solver of $p^1$ (2d t2\_1) in 2-d & JacobiKAN & FBiCGStab & 3-level & SSOR &  9-1-9 \\ \hline
Top 3 solver of $p^1$ (2d t3\_1) in 2-d & JacobiKAN & FBiCGStab & 3-level & SSOR & 7-1-7 \\ \hline
%Accurate-4($\frac{1}{128}$)& ChebyKAN & SSOR & 1-level & 0 & 0.015625 \\ \hline
%Accurate-5($\frac{1}{128}$)& DeepONet & BiCGStab & 1-level & 8 & 0.5 \\  \hline
\end{tabular}

    \end{minipage}

 %   \vspace{1em} % Space between subtables

    % Second subtable
    \begin{minipage}{\linewidth}
        \centering
        \captionof{subtable}{Performance of the top-3 solvers.}
         % Resize to fit the line width
          \begin{tabular}{c|c|c|c|c|c|c|c}
\hline \hline 
 &Error & Com. time &  \# of ite. & Conv. rate & Memory & MACs & Training time \\ 
\hline \hline
2d t1\_1& 1.20E-14 & 0.124 & 3 & 12.792 & 6826.64 & 2.131E+9 & 2207.45 \\ \hline
2d t2\_1&  2.58E-13 & 0.0758 & 2 & 13.873 & 7150.66 & 2.042E+8 & 3203.27  \\ \hline
2d t3\_1&  9.45E-15 & 0.1047 & 3 & 13.0879 & 7151.91 & 2.632E+8 & 3203.27 \\ \hline
%Accurate-4($\frac{1}{128}$)& 0.0971029 & 1.24E-14 & 6 & 2.86298 & 235.11 & 178858 & 99.2164 \\ \hline
%Accurate-5($\frac{1}{128}$)& 0.0518639 & 1.32E-14 & 13 & 1.47149 & 236.639 & 178858 & 326.378 \\  \hline
\end{tabular}
        
    \end{minipage}
    \label{top3_p1_2d_kry}
\end{table}

\begin{table}
\tiny
    \centering
    \caption{Top-3 solvers by preference function $p^2$ in 1-D for Krylov-based methods.}
    \begin{minipage}{\linewidth}
        \centering
        \captionof{subtable}{The top-3 solvers.}
       
        \begin{tabular}{c|c|c|c|c|c}
\hline \hline 
 &Neural operator & Classical solver & Multi-grid & Relaxation & Strategies \\ 
\hline \hline
Top 1 solver of $p^2$ (1d t1\_2) in 1-d & JacobiKAN & FBiCGStab & 1-level & SSOR & 9-1-9 \\ \hline
Top 2 solver of $p^2$ (1d t2\_2) in 1-d & JacobiKAN & FBiCGStab & 1-level & SSOR &  7-1-7 \\ \hline
Top 3 solver of $p^2$ (1d t3\_2) in 1-d & DeepONet & FBiCGStab & 1-level & SSOR & 7-1-7 \\ \hline
%Accurate-4($\frac{1}{128}$)& ChebyKAN & SSOR & 1-level & 0 & 0.015625 \\ \hline
%Accurate-5($\frac{1}{128}$)& DeepONet & BiCGStab & 1-level & 8 & 0.5 \\  \hline
\end{tabular}

    \end{minipage}

  %  \vspace{1em} % Space between subtables

    % Second subtable
    \begin{minipage}{\linewidth}
        \centering
        \captionof{subtable}{Performance of the top-3 solvers.}
         % Resize to fit the line width
          \begin{tabular}{c|c|c|c|c|c|c|c}
\hline \hline 
 &Error & Com. time &  \# of ite. & Conv. rate & Memory & MACs & Training time \\ 
\hline \hline
1d t1\_2& 1.95E-14 & 7.441E-4 & 2 & 15.181 & 44.972 & 200268 & 56.618 \\ \hline
1d t2\_2& 2.84E-14 & 7.748E-4 & 2 & 19.189 & 43.703 & 218748 & 56.618 \\ \hline
1d t3\_2& 1.50E-14 & 7.456E-4 & 2 & 20.046 & 43.029 & 3964910 & 174.409 \\ \hline
%Accurate-4($\frac{1}{128}$)& 0.0971029 & 1.24E-14 & 6 & 2.86298 & 235.11 & 178858 & 99.2164 \\ \hline
%Accurate-5($\frac{1}{128}$)& 0.0518639 & 1.32E-14 & 13 & 1.47149 & 236.639 & 178858 & 326.378 \\  \hline
\end{tabular}

    \end{minipage}
    \label{top3_p2_1d_kry}
\end{table}

\begin{table}
\tiny
    \centering
    \caption{Top-3 solvers by preference function $p^2$ in 2-D for Krylov-based methods.}
    
    \begin{minipage}{\linewidth}
        \centering
        \captionof{subtable}{The top-3 solvers.}
        
       \begin{tabular}{c|c|c|c|c|c}
\hline \hline 
 &Neural operator & Classical solver & Multi-grid & Relaxation & Strategies \\ 
\hline \hline
Top 1 solver of $p^2$ (2d t1\_2) in 2-d & JacobiKAN & FBiCGStab & 3-level & SOR & 9-1-9 \\ \hline
Top 2 solver of $p^2$ (2d t2\_2) in 2-d & JacobiKAN & FBiCGStab & 3-level & SSOR &  3-1-3 \\ \hline
Top 3 solver of $p^2$ (2d t3\_2) in 2-d & JacobiKAN & FBiCGStab & 3-level & SOR & 3-1-3 \\ \hline
%Accurate-4($\frac{1}{128}$)& ChebyKAN & SSOR & 1-level & 0 & 0.015625 \\ \hline
%Accurate-5($\frac{1}{128}$)& DeepONet & BiCGStab & 1-level & 8 & 0.5 \\  \hline
\end{tabular}

    \end{minipage}

   % \vspace{1em} % Space between subtablesapp_numerics

    % Second subtable
    \begin{minipage}{\linewidth}
        \centering
        \captionof{subtable}{Performance of the top-3 solvers.}
       % Resize to fit the line width
          \begin{tabular}{c|c|c|c|c|c|c|c}
\hline \hline 
 &Error & Com. time &  \# of ite. & Conv. rate & Memory & MACs & Training time \\ 
\hline \hline
2d t1\_2& 1.14E-14 & 0.0747 & 3 & 11.947 & 7152.57 & 2.630E+8 & 3203.27\\ \hline
2d t2\_2& 9.50E-15 & 0.0656 & 3 & 10.663 & 7153.36 & 1.867E+8 & 3203.27 \\ \hline
2d t3\_2& 1.13E-14 & 0.0490 & 4 & 8.685 & 7155.27 & 2.305E+8 & 3203.27 \\ \hline
%Accurate-4($\frac{1}{128}$)& 0.0971029 & 1.24E-14 & 6 & 2.86298 & 235.11 & 178858 & 99.2164 \\ \hline
%Accurate-5($\frac{1}{128}$)& 0.0518639 & 1.32E-14 & 13 & 1.47149 & 236.639 & 178858 & 326.378 \\  \hline
\end{tabular}
        
    \end{minipage}
    \label{top3_p2_2d_kry}
\end{table}

\begin{figure}
\centering
\subfigure[1-d by $p^1$]{
\includegraphics[width=0.23\textwidth]{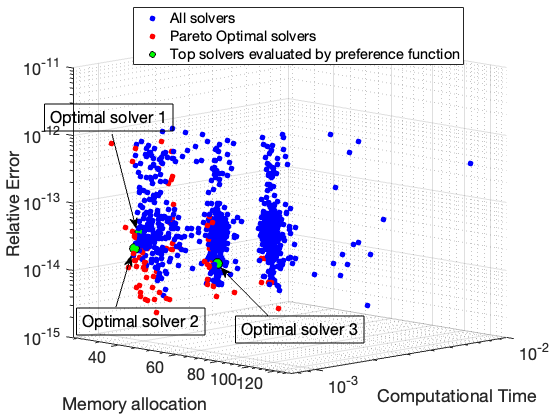}}
\subfigure[2-d by $p^1$]{
\includegraphics[width=0.23\textwidth]{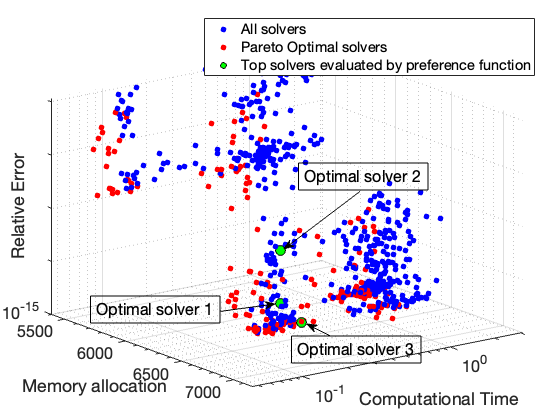}}
\subfigure[1-d by $p^2$]{
\includegraphics[width=0.23\textwidth]{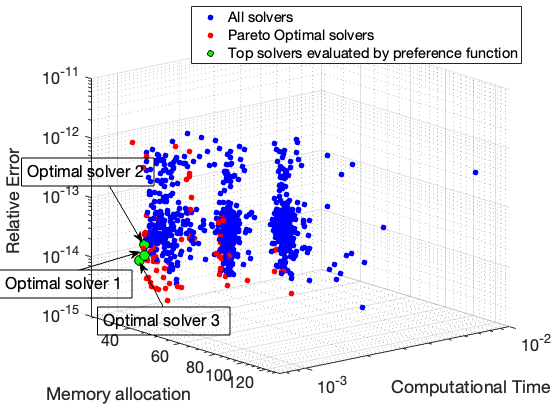}}
\subfigure[2-d by $p^2$]{
\includegraphics[width=0.23\textwidth]{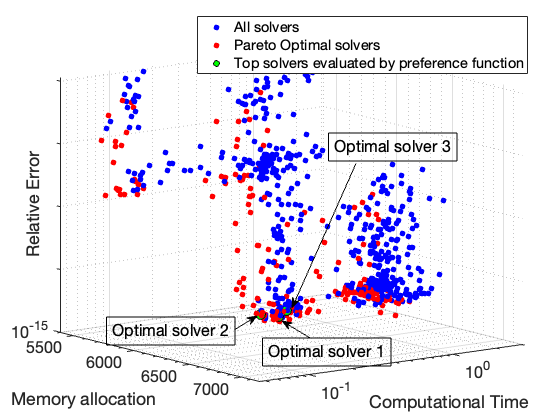}}

%\subfigure[Solving 3d Poisson Equation]{
%\includegraphics[width=0.3\textwidth]{figs/p1_kry_3d.png}}
\caption{The top-3 optimal solvers discovered by preference function $p^1$ (left 2) and $p^2$ (right 2), in the three-dimensional projection: Computational Time -- Relative Error -- Memory allocation of Pareto front, for Krylov-based methods.}
\label{front_p1_kry}
\end{figure}

\eqref{redis_kry_1d} represents the rediscovery of optimal solvers in 1-d case for Krylov-based methods, and \eqref{redis_kry_2d} represents the rediscovery of optimal solvers in 2-d.
\begin{subequations}\label{redis_kry_1d}
\begin{equation}
\left\{
\begin{aligned}
&\text{Preference function $p(\lambda_a;r) = (0.63, 0, 0, 0.12, 0.19, 0,  0.06)^Tr$ } \\
&\text{Optimal solver}: x^a = (\text{DeepONet}, \text{FGMRES}, \text{1-level}, \text{SOR}, 1-1-1) \ .
\end{aligned}
\right.
\end{equation}

\begin{equation}
\left\{
\begin{aligned}
&\text{Preference function $p(\lambda_a;r) = (0, 0, 0, 0.44, 0.56, 0,  0)^Tr$ } \\
&\text{Optimal solver}: x^a = (\text{DeepONet}, \text{FBiCGStab}, \text{1-level}, \text{SOR}, 1-1-1) \ .
\end{aligned}
\right.
\end{equation}

\end{subequations}

The following represent the rediscovery of optimal solvers in 2-d case:

\begin{subequations}\label{redis_kry_2d}

\begin{equation}
\left\{
\begin{aligned}
&\text{Preference function $p(\lambda_a;r) = (0.02, 0.51, 0, 0.23, 0.24, 0, 0)^Tr$ } \\
&\text{Optimal solver}: x^a = (\text{DeepONet}, \text{FCG}, \text{1-level}, \text{SSOR}, 7-1-7) \ .
\end{aligned}
\right.
\end{equation}

\begin{equation}
\left\{
\begin{aligned}
&\text{Preference function $p(\lambda_a;r) = (0.08, 0.04, 0, 0.35, 0.41, 0.11, 0.02)^Tr$ } \\
&\text{Optimal solver}: x^a = (\text{JacobiKAN}, \text{FGMRES}, \text{1-level}, \text{SSOR}, 9-1-9) \ .
\end{aligned}
\right.
\end{equation}

\end{subequations}